\documentclass[11pt]{amsart}
\usepackage{amssymb,amscd,amsthm,amsmath,amsfonts,epsf,latexsym}
\usepackage[dvips]{graphicx}

\input epsf

\newtheorem{thm}{Theorem}[section]
\newtheorem{theorem}[thm]{Theorem}
\newtheorem{cor}[thm]{Corollary}

\newtheorem{lemma}[thm]{Lemma}

\newtheorem{prop}[thm]{Proposition}

\theoremstyle{remark}
\newtheorem{remark}[thm]{Remark}

\theoremstyle{definition}
\newtheorem{defn}[thm]{Definition}

\numberwithin{equation}{section}

\newcommand{\g}{\mathfrak{g}}

\newcommand{\R}{{\mathbb{R}}}
\newcommand{\C}{{\mathbb{C}}}
\newcommand{\Q}{{\mathbb{Q}}}

\newcommand{\Z}{{\mathbb{Z}}}
\newcommand{\G}{{\mathbb{G}}}
\renewcommand{\O}{{\mathcal O}}
\newcommand{\X}{{\mathcal X}}

\newcommand{\GL}{\mathrm{GL}}
\newcommand{\Gal}{\mathrm{Gal}}
\newcommand{\End}{\mathrm{End}}
\newcommand{\Aut}{\mathrm{Aut}}
\newcommand{\Hom}{\mathrm{Hom}}
\newcommand{\SU}{\mathrm{SU}}
\newcommand{\SO}{\mathrm{SO}}
\newcommand{\SL}{\mathrm{SL}}
\newcommand{\PSL}{\mathrm{PSL}}
\newcommand{\Sp}{\mathrm{Sp}}
\newcommand{\Spin}{\mathrm{Spin}}
\newcommand{\diag}{\mathrm{diag}}

\newcommand{\1}{{1 \hspace{-0.35em} {\rm 1}}}
\newcommand{\I}{\mathcal{I}}
\newcommand{\A}{\mathcal{A}}

\newcommand{\Vn}{V^{\otimes n}}
\newcommand{\B}{\mathcal{B}}
\newcommand{\CC}{\mathcal{C}}
\newcommand{\M}{\mathcal{M}}

\newcommand{\lan}{\langle}
\newcommand{\ra}{\rangle}

\newcommand{\F}{\mathcal{F}}

\newcommand{\one}{\mathbf 1}

\newcommand{\N}{\mathbb{N}}

\newcommand{\la}{{\lambda}}

\begin{document}
\title[The $N$-eigenvalue Problem]
{The $N$-eigenvalue Problem and Two Applications}

\author{Michael J. Larsen}
\email{larsen@math.indiana.edu}
\address{Department of Mathematics\\
    Indiana University \\
    Bloomington, IN 47405\\
    U.S.A.}

\author{Eric C. Rowell}
\email{errowell@indiana.edu}
\address{Department of Mathematics\\
    Indiana University \\
    Bloomington, IN 47405\\
    U.S.A.}

\author{Zhenghan Wang}
\email{zhewang@indiana.edu}
\address{Department of Mathematics\\
    Indiana University \\
    Bloomington, IN 47405\\
    U.S.A.}

\thanks{The authors are partially supported by NSF grant  DMS-034772, and the third-named
author is also supported by NSF grant EIA 0130388.}

\begin{abstract}

We consider the classification problem for compact Lie groups
$G\subset U(n)$ which are generated by a single conjugacy class with
a fixed number $N$ of distinct eigenvalues.  We give an explicit
classification when $N=3$, and apply this to extract information
about Galois representations and braid group representations.
\end{abstract}
\maketitle

\section{Introduction}
\label{s:intro}

Consider the following two questions:
\begin{enumerate}
\item When does a compact Lie group
$G\subset U(n)$ have an element $g\in G$ possessing exactly two
eigenvalues.
\item When does a compact Lie group $G\subset U(n)$ have a cocharacter $U(1)\to G$ such that
the composition $U(1)\to U(n)$ is a representation of $U(1)$ with
exactly two weights.
\end{enumerate}

A solution to the second problem gives a family of solutions
to the first, by choosing $g$ to be almost any element of
the image of $U(1)$.  The converse is not true.  For one thing, any
non-central element of order $2$ in $G$ has exactly two eigenvalues.
To eliminate these essentially trivial solutions, we can  insist
that the ratio between the two eigenvalues is not $-1$.  There
remain interesting cases of finite groups $G$ satisfying the first
(but obviously not the second) condition, especially when the ratio
of eigenvalues is a third or fourth root of unity (see \cite{Bl},
\cite{Ko}, and \cite{Wa} for classification results).  On the other
hand, when $G$ is infinite modulo center, the solutions of the two
problems are essentially the same, though the historical reasons for
considering them were quite different. The first problem was
recently solved in the infinite-mod-center case by M.~Freedman,
M.~Larsen, and Z.~Wang \cite{FLW} with an eye toward understanding
representations of Hecke algebras.  The second problem was solved by
J-P. Serre \cite{Ser} nearly thirty years ago in order to classify
representations arising from Hodge-Tate modules of weight 1.

This paper is primarily devoted to an effort to understand the
analogue of the first problem (the ``$N$-eigenvalue problem'' of the
title) when the number $N\ge 3$ of eigenvalues is fixed and $G$ is
infinite modulo its center.  As a consequence, we also say something
about the second problem.  We are especially interested in the case
$N=3$, both because the results can be made quite explicit and
because it is especially relevant to the applications we have in
mind.  To specify our problem more precisely, we make the following
definitions.

A \emph{pair} $(G,V)$ consists of a compact Lie group $G$ and a
faithful irreducible complex representation $\rho\colon G\to\GL(V)$.
Let $N$ be a positive integer. We say a pair $(G,V)$ satisfies the
\emph{$N$-eigenvalue property} if there exists
a \emph{generating element}, i.e., an element $g\in G$ such that
the conjugacy class of $g$ generates $G$ topologically and the
spectrum $X$ of $\rho(g)$ has $N$ elements and satisfies the
\emph{no-cycle property}: for all roots of unity $\zeta_n$, $n\ge
2$, and all $u\in\C^\times$,
\begin{equation}
\label{e:cosets} u\langle \zeta_n\rangle \not\subset X
\end{equation}
Our goal is to classify pairs satisfying the $N$-eigenvalue property.

From the perspective of \cite{FLW}, the most obvious reason to
consider the $N$-eigenvalue property is that certain naturally
occurring representations of the Artin braid groups $\B_n$ satisfy
this condition.  The braid generators (half-twists) in the braid
group form a generating conjugacy class, and given any braided
tensor category $\CC$ and any object $x\in\CC$, we get a
representation of $\rho_{n, x}\colon \B_n\to \GL(V_{n,x})$. When
$\rho_{n, x}$ is unitary with respect to a hermitian form on
$V_{n,x}$, the closure of $\rho_{n, x}$ is a compact Lie group
endowed with a natural faithful representation and a generating
conjugacy class. It is often possible to control the eigenvalues of
half-twists, to guarantee the $N$-eigenvalue condition, and to
guarantee irreducibility. In the case when the braided tensor
category $\CC$ is modular, we obtain in addition representations of
the mapping class groups $\M(\Sigma_g)$ of closed oriented surfaces
$\Sigma_g$ for each genus $g$.  It is well-known that $\M(\Sigma_g)$
is generated by the (mutually conjugate) Dehn twists $D_c$ on $3g-1$
non-separating simple closed curves $c$ on $\Sigma_g$ (see
\cite{I}). If $\CC$ has $m$ simple object types, then each $D_c$ has
at most $m$ distinct eigenvalues as the eigenvalues of $D_c$ consist
of twists $\theta_i$ of the simple objects.  When the values
$\theta_i$ satisfy the no-cycle condition, it follows that each
irreducible constituent of the representation of $\M(\Sigma_g)$
arising from $\CC$ defines a pair satisfying the $N$-eigenvalue
property for some $N\leq m$.

The original motivation for the work of \cite{FLW} was for applications
to quantum computing.  In \cite{FKLW}, topological models of quantum
computing based on unitary topological quantum field theories
(TQFTs) are proposed. Given a topological model of quantum
computing, an important issue is whether or not this topological
model is capable of simulating the universal circuit model of
quantum computing \cite{NC}. This question actually depends on the
specific design of the topological quantum computer.  But for the
models based on braiding anyons in \cite{FKLW}, the
universality question is translated into a question about the
closures of the braid group representations. Quantum computing is
the processing of information encoded in quantum state vectors in
certain Hilbert spaces $V_{n}$ by unitary transformations.
Universality is the ability to efficiently move any state vector
$v\in V_n$ sufficiently close to any other state vector in $V_n$. A
theorem of Kitaev-Solovay (see \cite{NC}) guarantees efficiency if the
available unitary transformations in $U(V_n)$ form a dense subset of
$SU(V_n)$. Therefore, universality of topological models in \cite{FKLW}
is equivalent to the density of braid group representations.

The unitary Witten-Reshetikhin-Turaev Chern-Simons TQFTs based on
the gauge groups $\SU(N)$ and $\SO(N)$ are of particular interests due
to their relevance to braid statistics in condensed matter physics.
It was discovered in the 1980s that in dimension 2, there are
quasi-particles which are neither fermions nor bosons
\cite{Wi}. The most interesting of these \emph{anyons} are
non-abelian: when two such quasi-particles are exchanged,
their wave function is changed by a unitary matrix, rather than a
complex number, which depends on the exchanging paths (braiding). It
is predicted by physicists that the braid statistics of
quasi-particles in certain fractional quantum Hall liquids are
described by Jones' unitary braid groups representations or
equivalently the braid representations coming from the $\SU(2)$
TQFTs. Physicists have also proposed models of braid statistics
based on the $\SO(3)$ \cite{FF} and $\SO(5)$ \cite{Wn} TQFTs. Therefore, it may well be the
case that both the Jones and the BMW braid group representations
describe braid statistics of quasi-particles in nature.
Experiments are proposed to confirm those predictions \cite{DFN}.

Problem (2) is significant partly because of its relation to problem (1), but
in addition, there are number-theoretic applications, in the spirit of
\cite{Ser}.  We mention a global one: assuming the Fontaine-Mazur
conjecture, we can prove that if $K$ is a number field, $\bar K$ is
an algebraic closure of $K$, $G_K = \Gal(\bar K/K)$, and $X$ is a
non-singular projective variety over $K$, then $E_8$ does not occur
as a factor of the identity component of the Zariski-closure of
$G_K$ in the second \'etale cohomology group of $\bar X$.

The paper is organized as follows.  The second and third sections
treat the infinite-mod-center case of the $N$-eigenvalue problem.
The second section gives the general shape of the solution for all
$N$, and the third section gives an actual list for $N=3$.  The
fourth section shows that a fairly weak hypothesis on the actual
eigenvalues is enough to guarantee that $G$ is infinite modulo its
center.  The fifth section gives applications to number theory, and
the sixth section gives applications to braid group representations.
We conclude with a discussion of future applications to topology and
quantum computing.

\subsection*{Acknowledgements} The second-named
author would like to thank Hans Wenzl for many helpful
correspondences.
\section{Infinite groups}
\label{s:inf} In this section, we consider the general
$N$-eigenvalue problem for infinite compact groups. Our methods come
directly from \cite{FLW} and \cite{LW}.

\begin{lemma}
\label{l:permute} Let $V = V_1\oplus \cdots\oplus V_k$ be a complex
vector space and $T\colon V\to V$ a linear transformation permuting
the summands $V_i$ non-trivially. Then the spectrum of $V$ does not
satisfy (\ref{e:cosets}).
\end{lemma}

\begin{proof}
Renumbering if necessary, we may assume $V$ permutes
$V_1,V_2,\ldots,V_r$ cyclically, where $r\ge 2$.  Let $W =
V_1\oplus\cdots\oplus V_r$, let $\zeta_r = e^\frac{2\pi i}r$, and
let $S\colon W\to W$ act as the scalar $\zeta_r^i$ on $V_i$.  Then
$$S T\vert_W S^{-1} = \zeta_p T|_W,$$
so the spectrum of $T\vert_W$ is invariant under multiplication by
$\zeta_p$.  It is therefore a union of $\langle \zeta_p\rangle$-cosets.
\end{proof}

\begin{lemma}\label{decomplemma}
Let $(G_1,V_1)$ and $(G_2,V_2)$ be pairs and let $G$ denote the
image of $G_1\times G_2$ in $\GL(V_1\otimes V_2)$.  If $G$ satisfies
the $N$-eigenvalue property, then there exist integers $N_1$ and
$N_2$ such that $N_1+N_2-1\le N$, and subgroups $G'_1<G_1$ and
$G'_2<G_2$ such that $(G'_i,V_i)$ satisfies the $N_i$-eigenvalue
property  and $G'_i Z(G_i) = G_i$ for $i=1,2$.
\end{lemma}

\begin{proof}
Let $g\in G$ be a generating element, and let $(g_1,g_2)\in
G_1\times G_2$ map to $G$. Let $G'_i$ denote the subgroup of $G_i$
generated by the conjugacy class of $g_i$. As the conjugacy class of
$g$ generates $G$, $(\ker G_i\to G) G'_i = G_i$. By construction,
$\ker G_i\to G\subset Z(G_i)$.

The spectrum of $\rho(g)$ is the product of the spectra of
$\rho_i(g_i)$. So the lemma reduces to the following claim: if $X_1$
and $X_2$ are finite subgroups of an abelian group $A$ such that
$X_1+X_2$ does not contain a coset of a non-trivial subgroup of $A$,
then $|X_1+X_2|\ge |X_1|+|X_2|-1$. This is well-known (see, e.g.,
\cite{Ke}).
\end{proof}

If $(G,V)$ arises in this way, we say it is \emph{decomposable};
otherwise, it is \emph{indecomposable}.  Note that the tensor
product of pairs which satisfy the $N_1$ and $N_2$-eigenvalue
conditions need not satisfy the $N_1+N_2-1$-eigenvalue condition.
For one thing, the product of sets of cardinality $N_1$ and $N_2$
could be as large as $N_1 N_2$.  For another, the product of sets
satisfying the no-cycle property may itself fail to satisfy the
no-cycle property.

\begin{prop}
\label{indecomp} Let $(G,V)$ be an indecomposable pair.  If $G$ is
infinite modulo its center and $(G,V)$ satisfies the $N$-eigenvalue property for some
$N$, then $G = G^\circ Z(G)$.
\end{prop}

\begin{proof}
Let $g\in G$ be a generating element.  Then $g\in G^\circ$ implies
$G = G^\circ$, in which case there is nothing to prove.  If
$V|_{G^\circ}$ is not isotypic, then $g$ acts non-trivially on the
isotypic factors, and by Lemma~\ref{l:permute}, the spectrum of $g$
fails to satisfy property~(\ref{e:cosets}).  If $V|_{G^\circ} = W^n
= W\otimes U$, where $G^\circ$ acts trivially on $U$ and irreducibly
on $W$, then the span of $\rho(G^\circ)$ is $\End(W)\otimes
\mathrm{Id}_U\subset \End(V)$, so $\rho(G)$ lies in the normalizer
of $\End(W)\otimes \mathrm{Id}_U$, which is $\End(W)\End(U)$.  Thus
$\rho$ maps $G$ to $(\GL(W)\times\GL(U))/\C^\times$.  Let $\tilde G$
denote the cartesian square
\begin{equation*}
\begin{CD}
\tilde G @>\tilde\rho>> \GL(W)\times\GL(U) \\
@V\pi VV                       @VVV\\
G     @>\rho>> (\GL(W)\times\GL(U))/\C^\times.\\
\end{CD}
\end{equation*}
If $\tilde g\in\pi^{-1}(g)$, then the projections of $\tilde
\rho(\tilde g)$ to $\GL(W)$ and $\GL(U)$ have spectra satisfying the
no-cycle property, since the product of these spectra is the
spectrum of $\tilde \rho(\tilde g)$.  If $\dim W$ and $\dim U$ are
both $\ge 2$, then $(G,V)$ is decomposable, contrary to hypothesis.
As $G^\circ$ is not in the center of $G$, $\dim W\ge 2$.
It follows that $\dim U = 1$, i.e., the restriction of $V$ to
$G^\circ$ is irreducible.   Thus every element of $G$ which commutes with
$G^\circ$ lies in $Z(G)$.

It follows that for every $g\not\in G^\circ Z(G)$,
conjugation by $g$ induces an automorphism of $G^\circ$ which is not
inner.  By \cite[7.5]{St}, this implies that there exists a maximal
torus $T$ of $G^\circ$ such that $gTg^{-1} = T$ but conjugation by
$g$ induces a non-trivial automorphism of $T$.  The characters of
$T$ appearing in $V|_T$ span $X^*(T)\otimes\R$ since $V$ is a
faithful representation.  Therefore, a non-trivial automorphism of
$T$ must permute the weights of $V$ non-trivially.  By
Lemma~\ref{l:permute}, this implies that the spectrum of $g$
violates the no-cycle property, contrary to hypothesis.  Thus $g\in
G^\circ Z(G)$, and since the conjugacy class of $g$ generates $G$,
it follows that $G = G^\circ Z(G)$.
\end{proof}

\begin{prop}
\label{total-weight} Let $(G,V)$ be as in
Proposition~\ref{indecomp}.  Then $G$ is the product of the derived
group $D$ of $G^\circ$ and a group of scalar matrices in $V$.  The
group $D$ is simple modulo its center, and the restriction of $V$ to
$D$ is irreducible.  If the highest weight $\lambda$ of $V|_D$ is
written as a linear combination $\sum_i a_i\varpi_i$, where
$\varpi_i$ are the fundamental weights, then $\sum_i a_ib_i\le N-1$,
where the $b_i$ are positive integers determined by the root system of $D$.
\end{prop}

\begin{proof}
As $G^\circ$ is connected, $G^\circ = D Z(G^\circ)$.  As
$V|_{G^\circ}$ is irreducible, $Z(G^\circ)$ contains only scalars,
as does $Z(G)$.  Thus $G = DZ(G^\circ) Z(G)$, and the product
$Z(G^\circ) Z(G)$ is scalar in $\GL(V)$.  The centralizer of $D$ in
$\GL(V)$ equals the centralizer of $G^\circ = D Z(G^\circ)$ since
$Z(G^\circ)$ is scalar.  It follows that $V|_D$ is irreducible.
Any tensor decomposition of $V|_D$ extends to $G$ since scalars
respect any tensor decomposition; it follows that $V|_D$ is tensor
indecomposable and therefore that $D$ is simple modulo its center.
Let $\lambda$ denote its highest weight.

Let $g$ be a generating element, and let $t\in D$ be such that
$g^{-1}t$ is a scalar. $T$ be a maximal torus of $D$ containing $t$,
$R$ the root system of $D$ with respect to $T$, and $(\cdot,\cdot)$
the Killing form on $X^*(T)\otimes\R$.  Let
$$\langle\beta,\alpha\rangle=\frac{2(\beta,\alpha)}{(\alpha,\alpha)},$$
and fixing a Weyl chamber, let $\gamma$ denote the root dual to the
highest root in $R$. Thus $\gamma$ is the highest short root. By
\cite[VIII,~\S7,~Prop.~3(i)]{Bo}, the maximal arithmetic progression
of the form $\lambda,\lambda-\gamma,\lambda-2\gamma,\ldots$
contained in the set of weights of $V$ has length
$$1+\langle\lambda,\gamma\rangle = 1+\sum_i a_ib_i,$$
where the positive integers $b_i$ are the coefficients in the representation of
the highest root in $R$ in terms of the simple roots.
If this sum exceeds $N$, then the geometric progression of values
$$\lambda(t),\;(\lambda-\gamma)(t),\;(\lambda-2\gamma)(t),\;\ldots$$
must either take $\ge N+1$ distinct values, or fail
(\ref{e:cosets}), or be constant.  The first two possibilities are
ruled out by hypothesis, and it follows that $\gamma(t) = 1$.  If
$w$ belongs to the Weyl group, the same considerations apply to the
weight sequence $w(\lambda),w(\lambda)-w(\gamma),w(\lambda)-2
w(\gamma),\ldots$, so $w(\lambda)(t) = 1$.  On the other hand, the
short weights in a simple root system form a single Weyl orbit and
generate the root lattice, so $\alpha(t) = 1$ for all roots.  This
implies that $t$ lies in the center of $G$ and therefore that
$\rho(t)$ is scalar, contrary to hypothesis.
\end{proof}

One can also formulate the $N$-eigenvalue property for complex Lie
groups:

\begin{defn}
Let $G_{\C}$ be a reductive complex Lie group and $(\rho,V)$ a
faithful irreducible complex representation of $G_{\C}$. Then
$(G_{\C},V)$ satisfies the \emph{$N$-eigenvalue property} if there
exists a semisimple \emph{generating element} $g_{\C}\in G_{\C}$
whose conjugacy class generates a Zariski-dense subgroup
of $G_{\C}$, and such that the spectrum of
$\rho(g_{\C})$ consists of $N$ eigenvalues satisfying the no-cycle
condition.
\end{defn}

\begin{lemma}
Let $G_{\C}$ be a reductive complex Lie group and $(\rho,V)$ a
faithful irreducible complex representation of $G_{\C}$.  Let $G$ be
a maximal compact subgroup of $G_{\C}$.  Then $(G,V)$ satisfies the
$N$-eigenvalue property.
\end{lemma}

\begin{proof}
Let $T_{\C}$ denote the Zariski-closure of the cyclic group $\langle
g_{\C}\rangle$ and $T\subset T_{\C}$ the (unique) maximal compact
subgroup.  As $T$ can be regarded as the set of (real) points of a
real algebraic group whose complex points give $T_{\C}$, $T$ is
Zariski-dense in $T_{\C}$. We can decompose the restriction of $V$
to $T_{\C}$ as a direct sum of eigenspaces $V_\chi$ associated to
characters $\chi$ of $T_{\C}$.  There must be exactly $N$ such
eigenspaces, since any coincidence among
$\chi_1(g_{\C}),\ldots,\chi_{N+1}(g_{\C})$ gives the same
coincidence for the characters on all of $T_{\C}$.  The condition
that $\chi_i(t)\neq \chi_j(t)$ is open and non-empty in $T_{\C}$ as
is the condition that $\{\chi_1(t),\ldots,\chi_N(t)\}$ satisfy the
no-cycle condition.  It follows that $T$ contains an element $g$
which satisfies both conditions.

As all maximal compact subgroups of $G_{\C}$ are conjugate, without
loss of generality we may assume $T\subset G$.  We can regard $G$ as
the group of real points of a real linear algebraic group whose
complex points give $G_{\C}$ and $T\subset G$ as a Zariski-closed
subgroup. Let $H\subset G$ denote the smallest normal Zariski-closed
subgroup of $G$ containing $g$, or equivalently, $T$.  Thus $H$ can
be regarded as the group of real points of an algebraic group which
is a normal subgroup of the algebraic group with real locus $G$.
Let $H_{\C}$ denote the group of $\C$-points of this  subgroup.  If
$H\neq G$, then $H_{\C}\neq G_{\C}$, so $g_{\C}\in T_{\C}\subset
H_{\C}$ is contained in a proper normal subgroup of $g_{\C}$,
contrary to hypothesis. It follows that $g$ is a generating element
for $(G,V)$.

\end{proof}

\section{The $3$-eigenvalue problem}
\label{s:three}

In this section, we give an explicit solution of the
$\le3$-eigenvalue problem, assuming throughout that $G$ is a compact
Lie group which is infinite modulo center.

\begin{prop}
\label{p:two-eigenvalue} If $(G,V)$ is a pair satisfying the
2-eigenvalue property, and $\Phi$ denotes the root system of $G$ and
$\varpi$ the highest weight of $V$ in the notation of \cite{Bo},
then $(\Phi,\varpi)$ is one of the following:
\begin{enumerate}
\item $(A_r,\varpi_i)$, $1\le i\le r$.
\item $(B_r,\varpi_r)$.
\item $(C_r,\varpi_1)$.
\item $(D_r,\varpi_i)$, $i=1,r-1,r$.
\end{enumerate}
\end{prop}

\begin{proof}
This is the statement of \cite[1.1]{FLW}.
\end{proof}

Before treating the general $3$-eigenvalue problem, we make a
detailed study of the $A_r$ case.

\newcommand\cp{characteristic polynomial }
\begin{lemma}
\label{l:A-n} Let $(\rho,V)$ be an irreducible representation of
$\SU(n)$ with highest weight $\varpi$, and $t$ a non-central
element of $\SU(n)$. Suppose there are at most three eigenvalues of
$\rho(t)$ and they satisfy the no-cycle property. Then one of the
following is true:
\begin{enumerate}

\item For $1\le i\le n-1$, $\varpi = \varpi_i$,
and $t$ has \cp $(x-\lambda)^{n-1}(x-\lambda^{1-n})$; the
eigenvalues of $\rho(t)$ are $\lambda^i$, $\lambda^{i-n}$.
\label{i:two-mult-one}

\item For $1\le i\le n-1$, $\varpi = \varpi_i$,
and $t$ has \cp $(x-\lambda_1)^{n-2}(x-\lambda_2)^2$; the
eigenvalues of $\rho(t)$ are $\lambda_1^i$,
$\lambda_1^{i-1}\lambda_2$, and $\lambda_1^{i-2}\lambda_2^2 =
\lambda_1^{i-n}$. \label{i:two-mult-two}

\item For $i\in\{1,2,n-2,n-1\}$, $\varpi = \varpi_i$, and $t$ has
eigenvalues $\lambda_1$ and $\lambda_2$; the spectrum of $\rho(t)$ is $\{\lambda_1,\lambda_2\}$,
$\{\lambda_1^2,\lambda_1\lambda_2,\lambda_2^2\}$,
$\{\lambda_1^{-2},\lambda_1^{-1}\lambda_2^{-1},\lambda_2^{-2}\}$, or
$\{\lambda_1^{-1},\lambda_2^{-1}\}$, if $i$ is $1$, $2$, $n-2$, or
$n-1$ respectively. \label{i:any-two}

\item For $1\le i\le n-1$, $\varpi = \varpi_i$,
and $t$ has \cp
$(x-\lambda^{n-2})(x-\lambda\mu)(x-\lambda\mu^{-1})$; the
eigenvalues of $\rho(t)$ are $\lambda_1^i$, $\lambda_1^i\mu$, and
$\lambda_1^i\mu^{-1}$. \label{i:n-2-one-one}

\item For $i=1$ or $i=n-1$, $\varpi = \varpi_i$, and $t$ has eigenvalues
$\lambda_1,\lambda_2,\lambda_3$; the eigenvalues of $\rho(t)$ are the $\lambda_j$ or
the $\lambda_j^{-1}$ if $i=1$ or $i=n-1$ respectively.
\label{i:any-three}

\item For $i=1$ or $i=n-1$, $\varpi = 2\varpi_i$, and $t$ has eigenvalues
$\lambda_1$ and $\lambda_2$, each of multiplicity at least $2$; the eigenvalues of $\rho(t)$
are $\{\lambda_1^2,\lambda_1\lambda_2,\lambda_2^2\}$ or
$\{\lambda_1^{-2},\lambda_1^{-1}\lambda_2^{-1},\lambda_2^{-2}\}$ if
$i$ is $1$ or $n-1$ respectively. \label{i:sym-squared}

\item The highest weight $\varpi$ is $\varpi_1+\varpi_{n-1}$, and $t$ has eigenvalues
$\lambda_1$ and $\lambda_2$, each of multiplicity at least $2$; the
eigenvalues of $\rho(t)$ are $\lambda_1/\lambda_2$, $1$, and
$\lambda_2/\lambda_1$. \label{i:adjoint}

\item For $1\le i\le j\le n-1$, $\varpi = \varpi_i+\varpi_j$,
and $t$ has \cp $(x-\lambda)^{n-1}(x-\lambda^{1-n})$; the
eigenvalues of $\rho(t)$ are $\lambda^{i+j}$, $\lambda^{i+j-n}$, and
$\lambda^{i+j-2n}$. \label{i:any-ij}

\end{enumerate}
In particular, only case (\ref{i:any-three}) can give three
eigenvalues not in geometric progression.
\end{lemma}

\begin{proof}
By Proposition~\ref{total-weight}, if $\rho(t)$ has $N\le 3$
eigenvalues, $\varpi$ is a sum of at most $N-1$ fundamental weights.
If $\varpi = \varpi_i$ and $t$ has  eigenvalues
$\lambda_1,\ldots,\lambda_n$, the eigenvalues of $\rho(t)$ are
$$\left\{\prod_{s\in S}\lambda_s\Bigm\vert S\subset \{1,\ldots,n\}, |S| = i\right\}.$$
Duality exchanges $\varpi_i$ and $\varpi_{n-i}$ so without loss of
generality we may assume $i\le n/2$.  If
$\lambda_1,\ldots,\lambda_4$ are all distinct, and $n\ge i+3$ (in
particular, this holds if $n\ge 5$), then
$$\{\lambda_j\lambda_5\lambda_6\cdots\lambda_{3+i}\mid 1\le j\le 4\}$$
already contains four distinct elements.  If $n=4$ and $i=2$, two
products $\lambda_i\lambda_j$ and $\lambda_k\lambda_l$ are distinct
unless $\{i,j\}$ and $\{k,l\}$ are complementary sets, in which case
the equality implies $\lambda_i\lambda_j = \pm 1$.  At least one of
$\lambda_1\lambda_j$, $2\le j\le 4$ is neither $1$ nor $-1$, so
there must be at least four elements in the set
$\{\lambda_1\lambda_2,\ldots,\lambda_3\lambda_4\}$. If
$$\lambda_1=\lambda_2=\lambda_3\neq \lambda_4=\lambda_5=\lambda_6,$$
and $i\ge 3$, then
$$\{\lambda_1^j\lambda_4^{3-j}\lambda_7\cdots\lambda_{3+i}\mid 0\le j\le 3\}$$
contains a non-constant 4-term geometric progression in the spectrum
of $\rho(t)$, contrary to hypothesis.
If
$$\lambda_1=\lambda_2\neq\lambda_3=\lambda_4\neq\lambda_5\neq\lambda_1,$$
then
\begin{equation*}
\begin{split}
\{\lambda_1^2\lambda_2\lambda_6\cdots\lambda_{2+i},
       \lambda_1\lambda_2^2\lambda_6\cdots\lambda_{2+i},
       &\lambda_1\lambda_2\lambda_3\lambda_6\cdots\lambda_{2+i}, \\
       &\lambda_1^2\lambda_3\lambda_6\cdots\lambda_{2+i},
       \lambda_2^2\lambda_3\lambda_6\cdots\lambda_{2+i}\} \\
\end{split}
\end{equation*}
contains at least four distinct elements unless $\lambda_1\lambda_3
= \lambda_2^2$ and $\lambda_2\lambda_3 = \lambda_1^2$, in which case
it does not satisfy (\ref{e:cosets}). The remaining possibilities
are that $t$ has two distinct eigenvalues, one of multiplicity 1;
two distinct eigenvalues, one of multiplicity 2; two distinct
eigenvalues of arbitrary multiplicity, and $i$ (or $n-i$) is $\le
2$; three distinct eigenvalues, two of them of multiplicity 1; or
three distinct eigenvalues of arbitrary multiplicity, and $i$ (or
$n-i$) is $1$.

These give rise to cases (\ref{i:two-mult-one}),
(\ref{i:two-mult-two}), (\ref{i:any-two}), (\ref{i:n-2-one-one}),
and (\ref{i:any-three}) respectively. If $\lambda=
\varpi_i+\varpi_j$, $i\le j$, is among the weights appearing in
$V_\varpi$, then $\varpi_{i-1}+\varpi_{j+1}$ also appears, where we
define $\varpi_0=\varpi_n=0$. Thus if $\varpi = \varpi_i +
\varpi_j$, $i\le j$, then either $\varpi_{i+j}$, $\varpi_{2n-i-j}$
or $\varpi_1+\varpi_{n-1}$ is among the weights of $V_\varpi$, as
$i+j$ is less than, greater than, or equal to $n$.

First we consider the case $i+j=n$. If $t$ has three distinct
eigenvalues $\lambda_1,\lambda_2,\lambda_3$, then
$$|\{\lambda_1/\lambda_2,\lambda_2/\lambda_1,\lambda_1/\lambda_3,\lambda_3/\lambda_1,\lambda_2/\lambda_3,\lambda_3/\lambda_2\}|\le 3$$
implies that the set violates (\ref{e:cosets}) with $n=3$.  Thus,
$t$ has exactly two eigenvalues $\lambda_1$ and $\lambda_2$.  If
$i\ge 2$ and $\lambda_1$ and $\lambda_2$ each occurs with
multiplicity $\ge 2$, then
$$\{\lambda_1^2/\lambda_2^2,\lambda_1/\lambda_2,1, \lambda_2/\lambda_1,\lambda_2^2/\lambda_1^2\}$$
is contained in the spectrum of $\rho(t)$ since the Weyl orbits of
$\varpi_1+\varpi_{n-1}$ and $\varpi_2+\varpi_{n-2}$ are subsets of
the weights of $V_\varpi$.  As $\lambda_1\neq\lambda_2$, either this
set contains $5$ distinct elements or it violates (\ref{e:cosets}).
The remaining cases are (\ref{i:adjoint}) and the $i+j=n$ case of
(\ref{i:any-ij}).

If $i+j\neq n$, replacing $V_\varpi$ by its dual if necessary, we
can assume that $i+j < n$. If $3\le i+j\le n-3$, then $\varpi_{i+j}$
is a weight of $V_\varpi$, so by the analysis above, $t$ has two
eigenvalues, one with multiplicity one, and we are in case
(\ref{i:any-ij}).  If $i+j=2$, we see that $2\varpi_1$ and
$\varpi_2$ are both weights of $V_\varpi$, so if
$\lambda_1,\lambda_2,\lambda_3$ are eigenvalues of $t$,
$$\{\lambda_1^2,\lambda_2^2,\lambda_3^2,
\lambda_2\lambda_3,\lambda_3\lambda_1,\lambda_1\lambda_2\}$$
is contained in the spectrum of $\rho(t)$, contrary to assumption.
If there are exactly two eigenvalues, we get (\ref{i:sym-squared})
and the $i=j=1$ case of (\ref{i:any-ij}).

If $i+j=n-2$, then $V_\varpi$ contains all the weights of
$V_{\varpi_{n-2}}$, so $t$ may have only two eigenvalues,
$\lambda_1$ and $\lambda_2$, by the analysis of the case that
$\varpi$ is a fundamental weight, above.  If each occurs with
multiplicity $\ge 2$ and (without loss of generality) $\lambda_1$
occurs with multiplicity $\ge 3$, then
$$\{\lambda_2/\lambda_1^3,1/\lambda_1^2,1/\lambda_1\lambda_2,1/\lambda_2^2\}$$
is a $4$-term geometric progression contained in the spectrum of
$\rho(t)$ contrary to hypothesis. If $i+j=n-1$, then $V_\varpi$
contains all the weights of $V_{\varpi_{n-1}}$. If $\lambda_1$ and
$\lambda_2$ are eigenvalues of $t$ of multiplicity $\ge 2$, then the
spectrum of $\rho(t)$ contains the 4-term geometric progression
$$\{\lambda_2/\lambda_1^2, 1/\lambda_1,1/\lambda_2,\lambda_1/\lambda_2^2\}.$$
If $t$ has three distinct eigenvalues
$\lambda_1,\lambda_2,\lambda_3$, then the spectrum of $\rho(t)$
contains
$$\{1/\lambda_1,1/\lambda_2,1/\lambda_3,\lambda_1/\lambda_2\lambda_3,
\lambda_2/\lambda_1\lambda_3, \lambda_3/\lambda_1\lambda_2\}$$
which either violates the no-cycle condition or contains more than
$3$ elements. It follows that $t$ has exactly two eigenvalues, one
of multiplicity $n-1$. So all of these possibilities are subsumed in
case (\ref{i:any-ij}).

\end{proof}

\begin{thm}\label{mainthm}
If $(G,V)$ is an indecomposable pair satisfying the $3$-eigenvalue
property, $\Phi$ denotes the root system of the derived group $D$ of
$G^\circ$, and $\varpi$ the highest weight of $V$, then
$(\Phi,\varpi)$ is either one of the pairs enumerated in
Proposition~\ref{p:two-eigenvalue} or one of the following:
\begin{enumerate}
\item $(A_r,\varpi_i+\varpi_j)$, $1\le i\le j\le r$.
\item $(B_r,\varpi_i)$, $1\le i \le r-1$.
\item $(B_r,2\varpi_r)$.
\item $(C_r,\varpi_i)$, $2\le i\le r$.
\item $(C_r,2\varpi_1)$.
\item $(D_r,\varpi_i)$, $2\le i\le r-2$.
\item $(D_r,\varpi)$, $\varpi\in\{2\varpi_{r-1},\varpi_{r-1}+\varpi_r,2\varpi_r\}$.
\item $(E_6,\varpi_i)$, $i=1,3,6$.
\item $(E_7,\varpi_i)$, $i=1,7$.
\item $(F_4,\varpi_4)$.
\item $(G_2,\varpi_2)$.
\end{enumerate}
If there exists a generating element with three eigenvalues which do
not form a geometric progression, then $(\Phi,\varpi)$ is
$(A_r,\varpi_1)$ or $(A_r,\varpi_r)$.
\end{thm}

\begin{proof}
By Proposition~\ref{total-weight}, the root system is simple and if
$\varpi = \sum_i a_i\varpi_i$ and the highest root is $\sum_i
b_i\alpha_i$, then $\sum a_i b_i \le 2$. By \cite[Planches]{Bo},
this reduces the possibilities to those listed, together with:
\begin{enumerate}
\setcounter{enumi}{11}
\item $(D_r,\varpi)$, $\varpi\in\{2\varpi_1,\varpi_1+\varpi_{r-1},\varpi_1+\varpi_r\}$.
\label{bad-D}
\item $(E_6,\varpi)$, $\varpi\in\{2\varpi_1, \varpi_2,\varpi_5,2\varpi_6,\varpi_1+\varpi_6\}$.
\label{bad-E6}
\item $(E_7,\varpi)$, $\varpi\in\{\varpi_2,\varpi_6,2\varpi_7\}$.
\label{bad-E7}
\item $(E_8,\varpi)$, $\varpi\in\{\varpi_1,\varpi_8\}$
\label{bad-E8}
\item $(F_4,\varpi_1)$.
\label{bad-F4}
\end{enumerate}

To see that the classical cases (1)--(7) above are achieved, we let
$G = D$ and $V$ the indicated representation, and we choose the
generating element as follows.  For $A_r$, we let $g$ be the image
of the diagonal element
$\diag(\lambda^{-r},\lambda,\cdots,\lambda)\in\SU(r+1)$ in $G$.  For
$B_r$, we let $g$ denote the image of an element in $\Spin(2r+1)$
whose image in $\SO(2r+1)$ is
$\diag(\lambda,\lambda^{-1},1,\ldots,1)$ For $C_r$, we let $g$
denote the image of the element $(\lambda,\lambda^{-1},1\ldots,1)$
in $\Sp(2r)$. For $D_r$, we let $g$ denote the image of an element
in $\Spin(2r)$ whose image in $\SO(2r)$ is
$\diag(\lambda,\lambda^{-1},1,\ldots,1)$.

Next we show that the excluded cases (\ref{bad-D})--(\ref{bad-F4}) above do not occur.
For $D_r$, we consider an element $g$ whose image in $\SO(2r)$ has
eigenvalues $\lambda_1^{\pm 1}$, $\ldots$, $\lambda_r^{\pm 1}$. In
$V_{2\varpi}$, the eigenvalues of $g$ are $\lambda_i^{\pm 2}$,
$\lambda_i^{\pm 1}\lambda_j^{\pm 1}$, and 1. It is easy to see these
represent at least 5 distinct values.  A similar analysis rules out
the remaining cases in (\ref{bad-D}).

For $E_6$ and $F_4$ we use the existence of equal rank semisimple
subgroups of the form $A_2^k$.  As these subgroups share a maximal
torus with their ambient groups, every generating element $g$ can be
conjugated into the subgroup. We use the branching rules tabulated
in \cite{MP} to compute the restrictions of $G$-representations via
$\SU(3)^k\to  G$; since the center of $\SU(3)^k$ has exponent $3$,
and since we know that there are no $2$-eigenvalue solutions for
$F_4$ and $E_6$, there can be no $\le 3$-eigenvalue solutions coming
from central elements of $\SU(3)^k$ and satisfying (\ref{e:cosets}).
If $M(\lambda)\in\SU(3)$ has eigenvalues
$\lambda,\lambda,\lambda^{-2}$, then
$M(\lambda)\times M(\lambda^{-1})$ maps to an element of $F_4$ which has eigenvalues
$\lambda^{-3},1,\lambda^3$ for $V_{\varpi_4}$.  The restriction of
$F_4$ to $\SU(3)^2$ is
$$V_{2\mu_2}\boxtimes V_{\mu_1}\oplus V_{2\mu_1}\boxtimes V_{\mu_2}\oplus
V_{\mu_1+\mu_2}\boxtimes V_0\oplus V_0\boxtimes V_{\mu_1+\mu_2};$$
the image of any element non-central in both factors has at least
four eigenvalues from the first summand; the image of any element
central in the second factor but not the first has at least four
eigenvalues from the first two summands; the image of any element
central in the first factor but not in the second has at least four
eigenvalues or the eigenvalues $\{1,e^{\pm 2\pi i/3}\}$ from the
first two summands.  For $(E_6,\varpi_1)$, the image of
$M(\lambda)\times M(\lambda)\times 1$ has eigenvalues
$\{\lambda^{-2},\lambda,\lambda^4\}$, and it is not difficult to see
that this is essentially the only way to get three eigenvalues.  For
$(E_6,\varpi_2)$, the image of $M(\lambda)\times M(\lambda)\times 1$
has eigenvalues $\{\lambda^{-3},1,\lambda^3\}$.  To see that the
excluded cases (\ref{bad-E6}) do not give solutions to the
$3$-eigenvalue problem, we note that
$$V_{\varpi_2}|_{A_2^3}
= V_{\mu_1+\mu_2}\boxtimes V_{\mu_2}\boxtimes V_{\mu_2} \oplus
V_{\mu_2}\boxtimes V_{\mu_1+\mu_2}\boxtimes V_{\mu_1}
\oplus\cdots;$$
$$V_{2\varpi_1}|_{A_2^3}
= V_{\mu_1+\mu_2}\boxtimes V_{\mu_2}\boxtimes V_{\mu_2} \oplus
V_{\mu_2}\boxtimes V_{\mu_1+\mu_2}\boxtimes V_{\mu_1}
\oplus\cdots;$$
$$V_{\varpi_1+\varpi_6}|_{A_2^3}
= V_{2\mu_1}\boxtimes V_{\mu_1}\boxtimes V_{\mu_2} \oplus
V_{\mu_1}\boxtimes V_{2\mu_1}\boxtimes V_{\mu_1} \oplus\cdots.$$
These summands are already enough to guarantee that if
$(E_6,\varpi_2)$, $(E_6,2\varpi_1)$, or $(E_6,\varpi_1+\varpi_6)$
satisfies the $3$-eigenvalue condition, any generating element in
$A_2^3$ must be central in two of the three factors and have
eigenvalues $\lambda,\lambda,\lambda^{-2}$, $\lambda^3\neq 1$, in
the third.  However, if $\omega^3 = 1$, neither
$$\{\lambda^{-3},1,\lambda^3,\omega\lambda,\omega\lambda^{-2}\}$$
nor
$$\{\lambda^2,\lambda^{-1},\lambda^{-4},\omega\lambda,\omega\lambda^{-2}\}$$
can have order $\le 3$ and satisfy the no-cycle property.

For $E_n$, $n\ge 7$, we use the equal rank subgroups $A_n$.  Again,
\cite{MP} gives the restriction of $V_\varpi$ to $\SU(n+1)$.  The
following table lists all irreducible components of these
restrictions for all possible $\varpi$.  It also specifies the
eigenvalues in $V_\varpi$ for the image of the scalar matrix $\zeta
I$ and the matrix $M(\lambda)$:
\smallskip
\begin{tabular}{|l|l|p{1.4in}|l|p{1.4in}|}\hline
$\Phi$  &$\varpi$       &$\{\mu_i\}$    &$\zeta I$ e-values
&$M(\lambda)$ e-values \\ \hline $E_7$   &$\varpi_1$ &$\mu_1+\mu_7$,
$\mu_4$
    &$\pm 1$        &$\lambda^{-8}$, $\lambda^{-4}$, $1$, $\lambda^4$, $\lambda^8$ \\
$E_7$   &$\varpi_2$ &$\mu_1+\mu_5$, $\mu_3+\mu_7$, $2\mu_1$,
$2\mu_7$
    &$\pm i$        &$\lambda^{-14}$, $\lambda^{-10}$, $\lambda^{-6}$, $\lambda^{-2}$,
                    $\lambda^2$, $\lambda^6$, $\lambda^{10}$, $\lambda^{14}$ \\
$E_7$   &$\varpi_6$ &$\mu_1+\mu_3$, $\mu_5+\mu_7$, $\mu_1+\mu_7$,
$\mu_2+\mu_6$
    &$\pm 1$        &$\lambda^{-12}$, $\lambda^{-8}$, $\lambda^{-4}$,
                 $1$, $\lambda^4$, $\lambda^8$, $\lambda^{12}$\\
$E_7$   &$\varpi_7$ &$\mu_2$, $\mu_6$
    &$\pm i$        &$\lambda^{-6}$, $\lambda^{-2}$, $\lambda^2$, $\lambda^6$ \\
$E_7$   &$2\varpi_7$    &$0$, $\mu_4$, $\mu_2+\mu_6$, $2\mu_2$,
$2\mu_6$
    &$\pm 1$        &$\lambda^{-12}$, $\lambda^{-8}$, $\lambda^{-4}$, $1$,
                    $\lambda^4$, $\lambda^8$, $\lambda^{12}$\\ \hline
$E_8$   &$\varpi_1$ &$\mu_3$, $\mu_6$, $\mu_1+\mu_8$
    &$1,e^{\pm 2\pi i/3}$   &$\lambda^{-9}$, $\lambda^{-6}$, $\lambda^{-3}$, $1$,
                            $\lambda^3$, $\lambda^6$, $\lambda^9$ \\
$E_8$   &$\varpi_8$ &$\mu_1+\mu_2$, $\mu_1+\mu_5$, $\mu_1+\mu_8$,
$\mu_2+\mu_7$,
                                $\mu_4+\mu_8$, $\mu_7+\mu_8$
    &$1,e^{\pm 2\pi i/3}$   &$\lambda^{-15}$, $\lambda^{-12}$, $\lambda^{-9}$, $\lambda^{-6}$,
                        $\lambda^{-3}$, $1$, $\lambda^3$, $\lambda^6$,
                        $\lambda^9$, $\lambda^{12}$, $\lambda^{15}$\\ \hline
\end{tabular}
\smallskip

It follows that neither scalar matrices nor matrices of the form
$M(\lambda)$ give rise to $3$-eigenvalue solutions. By
Lemma~\ref{l:A-n}, the only possible solutions to the $3$-eigenvalue
problem for $E_7$ and $E_8$ are the pairs $(E_7,\varpi_1)$,
$(E_7,\varpi_7)$, and $(E_8,\varpi_1)$.  For the first, an element
of $\SU(8)$ with eigenvalues $\lambda$, $\lambda$, $\lambda$,
$\lambda$, $\lambda$, $\lambda$, $\lambda^{-3}$, $\lambda^{-3}$ maps
to an element of $E_7$ with eigenvalues $\lambda^{-4}$, $1$,
$\lambda^4$.  For the second, an element of $\SU(8)$ with
eigenvalues $\lambda$, $\lambda$, $\lambda$, $\lambda$,
$\lambda^{-1}$, $\lambda^{-1}$, $\lambda^{-1}$, $\lambda^{-1}$ maps
to an element of $E_7$ with eigenvalues $\lambda^{-2}$, $1$,
$\lambda^2$.  For $(E_8,\varpi_1)$, the only possibility is an
element of $\SU(9)$ with $\lambda_1$ of multiplicity $7$ and
$\lambda_2$ of multiplicity $2$. This maps to an element of $E_8$
with two three-term geometric progressions of eigenvalues:
$\lambda_1^3$, $\lambda_1^2\lambda_2$, $\lambda_1\lambda_2^2$; and
$\lambda_1\lambda_2^{-1}$, $1$, $\lambda_1^{-1}\lambda_2$.  To have
three eigenvalues in all, we must have $\lambda_1^3 =
\lambda_1\lambda_2^{-1}$, which together with
$\lambda_1^7\lambda_2^2 = 1$ implies that the eigenvalues are all
equal, which we have already seen is not a possibility.

The case of $G_2$ is trivial.

When there are three eigenvalues not in geometric progression, the
representations cannot be self-dual, and if $\phi=A_r$, then
$\varpi\in\{\varpi_1,\varpi_r\}$ by Lemma~\ref{l:A-n}.  The only
remaining cases for which $V_\varpi$ is not self-dual are
$(D_r,V_{2\varpi_{r-1}})$ and its dual (when $r$ is odd) and
$(E_6,\varpi_1)$ and its dual.  In the first case, as $r$ is odd,
the Weyl orbit of $\varpi_1$ lies in the set of weights of both
$V_{2\varpi_{r-1}}$ and $V_{2\varpi_r}$.  The eigenvalues
contributed by these weights come in mutually inverse pairs; if
there are $\le 3$ but not three in geometric progression, then there
must be two: $\lambda$ and $\lambda^{-1}$, which are distinct from
one another.  Then the Weyl orbit of $\varpi_3$ also lies in the set
of weights of $V_\varpi$, so $\lambda^3,\lambda,\lambda^{-1},
\lambda^{-3}$ are all eigenvalues of $\rho(t)$, which is absurd. In
the second case, restricting from $E_6$ to $\SU(6)\times \SU(2)$, we
get
$$V_{\varpi_1}\boxtimes V_{\varpi_1}\oplus V_{\varpi_4}\boxtimes V_0.$$
The second summand contributes $\le 2$ eigenvalues or 3 eigenvalues
not in geometric progression, so an inverse image $(g_1,g_2)$ of the
generating element must be a scalar $\zeta$ in  $\SU(6)$ (and
therefore $\zeta^6=1$). The eigenvalues of $g$ in the first summand
are $\{\zeta\lambda,\zeta^4,\zeta\lambda^{-1}\}$ which are in
geometric progression, contrary to assumption.

\end{proof}

\section{The asymptotic $N$-eigenvalue condition}
\label{s:asymp}

In this section we consider what can be said when the eigenvalues of
a generating element are sufficiently general.  One hypothesis which
is strong enough for our purposes is that the eigenvalues are
distinct $r$th roots of unity where $r$ is a sufficiently large
prime.  We consider a somewhat more general condition.

\begin{prop}
\label{p:nbd} Let $T\cong U(1)^d$ be a torus and $U$ an open
neighborhood of the identity in $T$. There exists a finite set $S$
of characters $\chi\colon T\to U(1)$ and an integer $m$ such that if
$n$ is a positive integer and $t\in T$ an $n$-torsion point, at
least one of the following must be true:
\begin{enumerate}
\item There exists $\chi\in S$ such that $\chi(t)\neq 1$ has order $\le m$.
\item There exists an integer $k$ relatively prime to $n$ such that $t^k\in U$.
\end{enumerate}
\end{prop}

\begin{proof}
We use induction on dimension, the proposition being trivial in
dimension 0.

By Urysohn's lemma there exists a continuous function $f\colon
T\to[0,1]$ such that $f(x) = 0$ for $x\not\in U$ and $f(x)=1$ in
some neighborhood of the identity.  It is well-known (see, e.g.
\cite[VII~Th.~1.7]{SW}) that finite linear combinations of
characters are dense in the $L^\infty$ norm on the set of continuous
functions on $T$.  It follows that there exists a real-valued finite
sum $f(x) := \sum_{\chi\in S} a_\chi \chi(x)$ such that $f(x) < 0$
for all $x\in T\setminus U$ and $a_0 = \int f(x)\,dx > 0$.
Enlarging $S$ if necessary, we may assume without loss of generality
that if $n\chi\in S$ for some positive integer $n$, then $\chi\in
S$.

Suppose $\chi(t) = 1$ for some non-trivial character $\chi\in S$.
Let $\lambda\in S$ denote a primitive character in $S$ and $k$ a
positive integer such that $\chi = k\lambda$.  If $m$ is taken
greater than the value of $k$ associated with any character in $S$,
either (1) is satisfied or $\lambda(t) = 1$.  As $\lambda$ is
primitive, $\ker\lambda$ is a subtorus of $T$.  As there are only
finitely many subtori arising in this way, the proposition follows
by induction.

We may therefore assume that the order of $\chi(t)$ is greater than
$m$ for each $\chi\in S$.  We have
$$\sum_{\{k\in [0,n]\cap\Z\mid (k,n)=1\}} f(t^k) = a_0\phi(n)
    + \sum_{\chi\in S\setminus\{0\}}a_\chi \sum_{\{k\in [0,n]\cap\Z\mid (k,n)=1\}}\chi(t^k).$$
If $n_\chi$ is the order of $\chi(t)$, then
$$\sum_{\{k\in [0,n]\cap\Z\mid (k,n)=1\}}\chi(t^k)
= \frac{\phi(n)}{\phi(n_\chi)}\sum_{\{k\in [0,n_\chi]\cap\Z\mid
(k,n_\chi)=1\}}\chi(t^k) = \frac{\mu(n)\phi(n)}{\phi(n_\chi)}.$$
Choosing $m$ large enough that for all $n_\chi > m$,
$$\sum_{\chi\in S\setminus\{0\}} |a_\chi| \le \phi(n_\chi) a_0,$$
we conclude that
$$\sum_{\{k\in [0,n]\cap\Z\mid (k,n)=1\}} f(t^k) \ge 0$$
and therefore that $t^k\in U$ for some $k$ prime to $n$.
\end{proof}

\begin{theorem}
For every integer $N\ge 2$ there exists an integer $m$ such that if
$(G,V)$ satisfies the $N$-eigenvalue property with a generator $g$
with eigenvalues $\lambda_1,\ldots,\lambda_N$, and $G$ is finite modulo its center,
then the group $\langle \lambda_i \lambda_j^{-1}\rangle$ generated
by ratios of eigenvalues of $\rho(g)$ contains a non-trivial root of
unity of order less than $m$.
\end{theorem}

\begin{proof}
If $G$ is finite modulo its center and acts irreducibly on $V$, then either $G^\circ$
is trivial or it consists of all scalars of absolute value 1.  In the latter case, we can
replace $g$ by $\det(g)^{1/\dim(V)} g$  for any choice of root, and the resulting conjugacy
class still satisfies the $N$-eigenvalue property, generates a subgroup of $G\cap \SL(V)$
(which is finite), and determines the same group of eigenvalue ratios
$\langle \lambda_i \lambda_j^{-1}\rangle$.  Without loss of generality, therefore, we may assume
$G$ is finite.

Any automorphism of $\C$ determines an automorphism of the abstract
group $\GL_n(\C)$ for each $n$.  Consider the quotient $T =
U(1)^n/U(1)$ of the diagonal unitary matrices by the unitary scalar
matrices.  Let $U\subset T$ denote the image of $A^n$ in $T$, where
$A$ is the arc from $-\pi/6$ to $\pi/6$, and let $n$ be the order of
the group generated by the eigenvalues of $g$. We apply
Proposition~\ref{p:nbd} to obtain $m$ large enough that our
hypotheses imply the existence of a field automorphism $\sigma$ of
$\C$ such that all the eigenvalues of $\sigma(\rho(g))$ lie in an
arc of length $\le \pi/3$ on the unit circle.  By \cite[Theorem
8]{Bl}, this implies that the representation representation
$\sigma{\scriptstyle\circ}\rho$ is imprimitive.  As the conjugacy
class of $g$ generates $G$, the element $g$ itself must satisfy the
hypothesis of Lemma~\ref{l:permute}, and therefore the spectrum of
$\sigma(\rho(g))$ does not satisfy (\ref{e:cosets}). As this property
is stable under Galois action, the spectrum of $\rho(g)$ fails to
satisfy (\ref{e:cosets}), contrary to hypothesis.
\end{proof}

\begin{cor}
For every integer $N\ge 2$ there exists an integer $m$ such that if
$(G,V)$ satisfies the $N$-eigenvalue property with a generator $g$
of prime order $r$, then $r < m$ or $G$ is infinite modulo its center.
\end{cor}

We remark that it is probably possible to prove a stronger version
of this corollary, in which a good bound is given for $m$, using
\cite{Za} as a starting point.

\section{Application to Hodge-Tate theory}

Let $\bar\Q_p$ be an algebraic closure of $\Q_p$, and $\C_p$ denote
the completion of $\bar\Q_p$.  Let $K$ and $L$ be subfields of
$\bar\Q_p$ finite over $\Q_p$, and let $\Gamma_K := \Gal(\bar\Q_p/K)$.
Let $V_L\cong L^d$ be a finite-dimensional $L$-vector space and
$\rho_L\colon \Gamma_K\to\GL(V_L)$ a continuous representation.
Then $\Gamma_K$ acts on both factors of $V_{\C_p}:=V_L\otimes_L\C_p$.
The representation is said to be \emph{Hodge-Tate} if
$V_{\C_p}$ decomposes as a direct sum of factors $V_{i\C_p}$ such
that $\Gamma_K$ acts on $V_i$ through the $i$th tensor power of the
cyclotomic character.  If $X$ is a complete
non-singular variety over $K$ and $\bar X$ is obtained from $X$ by
extending scalars to $\bar\Q_p$, then $V_L:=H^k(\bar X,L)$ is
Hodge-Tate for all non-negative integers $k$, and the factors
$V_{i\C_p}$ are non-zero only if $0\le i\le k$ (\cite{Fa}).

Let $G_L$ denote the Zariski-closure of the image of
$\rho_L(\Gamma_K)$ in $\GL_d$. By the axiom of choice, any two
uncountable algebraically closed fields of characteristic zero whose
cardinalities are the same are isomorphic.  Therefore,
$\C\cong\C_p$, and extending scalars, we can view $G_{\C}$ as a
complex algebraic group.  Let $G$ denote a maximal compact subgroup
of $G_{\C}$.  The inclusion $G_{\C}\subset \GL(V_{\C})$ gives $G$ a
complex representation which we denote $(\rho,V)$.  If $\rho_L$ is
absolutely irreducible, then $V_{\C}$ is an irreducible
representation of $G_{\C}$ and therefore of $G$.

Although $G_L$ need not be connected, by passing to a finite
extension $K'$ of $K$ (i.e., replacing $\Gamma_K$ by a normal open
subgroup) we can replace $G_L$ by its identity component.
Therefore, in trying to understand what Lie algebras and Lie algebra
representations can arise from Hodge-Tate structures with specified
weights, without loss of generality we may assume that $G_L$ is
connected.

\begin{defn}
Let $G_{\C}$ be a connected reductive algebraic group over $\C$, and
$V$ a faithful complex representation of $G_{\C}$. We say that
$(G_{\C},V)$ is of \emph{$N$-eigenvalue type} if for every almost
simple normal subgroup $H_{\C}$ of $G_{\C}$ and every irreducible
factor $W$ of $V|_{H_{\C}}$, the image of $H_{\C}$ in $\GL(W)$
satisfies the $N_W$-eigenvalue property for some $N_W\le N$.
\end{defn}

\begin{lemma}
\label{l:complex} Let $G_{\C}$ be a connected reductive complex Lie
group and $(\rho,V)$ a faithful representation. Let $g_{i\C}\in
G_{\C}$ be semisimple elements generating a Zariski-dense subgroup of
$G_{\C}$, such that the
spectrum of $\rho(g_{i\C})$ has $N$ eigenvalues satisfying the no-cycle
condition. Then $(G_{\C},V)$ is of $N$-eigenvalue type.
\end{lemma}

\begin{proof}
Let $D_{\C}$ denote the derived group of $G_{\C}$. The universal
cover $\tilde D_{\C}$ factors into simply connected, almost simple
complex groups $G_{j\C}$.  Every irreducible factor $W$ of $V$
restricts to an irreducible representation of $\tilde D_{\C}$ which
decomposes as $W_1\otimes\cdots\otimes W_k$, where $W_j$ is an
irreducible representation of $G_{j\C}$.

Each $g_{i\C}$ in our generating set factors as $d_{i\C} z_{i\C}$,
where $z_{i\C}$ lies in the center of $G_{\C}$. We choose $\tilde
d_{i\C}\in \tilde D_{\C}$ lying over $d_{i\C}$, and let $g_{ij\C}$
denote the $G_{j\C}$ coordinate of $\tilde d_{i\C}$.  For each $j$,
there exists $W$ such that $W_j$ is non-trivial and $i$ such that
$g_{ij\C}$ does not lie in the center of $G_{j\C}$. As $g_{i\C}$ is
semisimple, the same is true of $d_{i\C}$ and therefore $\tilde
d_{i\C}$ and therefore $g_{ij\C}$.  Moreover, it has at most $N$
eigenvalues on $W_j$ and they satisfy (\ref{e:cosets}), since
if $S$ and $T$ are sets of complex numbers and the product set
satisfies (\ref{e:cosets}), then $|S|,|T|\le |ST|$, and $|S|$ and
$|T|$ satisfy (\ref{e:cosets}). As $g_{ij\C}$ is not in the center of
$G_{j\C}$, the conjugacy class of $\rho_j(g_{ij\C})$ generates a
non-central normal subgroup of the almost simple group $\rho_j(G_{j\C})$
and therefore generates the whole group.

\end{proof}

\begin{theorem}
\label{t:HT} If $V_L$ is an absolutely irreducible Hodge-Tate
representation of $G_K$ with $N$ distinct weights, then
$(G_{\C}^\circ,V)$ is of $N$-eigenvalue type.
\end{theorem}

\begin{proof}
The grading of $V_{\C}$ which assigns $V_{i\C}$ degree $i$ uniquely
determines a cocharacter $h\colon\G_m\to G_{\C}$ such that
$\rho{\scriptstyle\circ} h$ acts isotypically on $V_{i\C}$ by the
$i$th power character.  By \cite{Sen}, $G_L^\circ$ is the smallest
$L$-algebraic subgroup of $\GL_d$ which contains $h(\G_m)$. Thus
$\{h^\sigma(\G_m)\mid\sigma\in\Aut_L(\C)\}$ generates
$G_{\C}^\circ$.  If $u\in\C^\times$ is of infinite order, then any element
$g_{j\C}\in h^{\sigma_j}(u)$ (Zariski-topologically) generates
$h^{\sigma_j}(\G_m)$.  Together, the $g_{j\C}$ generate $G_{\C}^\circ$.
There are exactly $N$ distinct eigenvalues of $\rho(g_{j\C})$ and they satisfy the
no-cycle condition.
The theorem now follows from  Lemma~\ref{l:complex}.
\end{proof}

\begin{theorem}
Assume that the Fontaine-Mazur conjecture \cite[Conj.~5a]{FM} holds.
If $X$ is a complete non-singular variety over a number field $K$,
$k$ is a non-negative integer, $G_{\C}$ is the complexification of
the Zariski closure of $\Gal(\bar K/K)$ in $\Aut(H^k(\bar X,\Q_p))$,
and $V = \Aut(H^k(\bar X,\Q_p))\otimes_{\Q_p}\C$, then
$(G_{\C}^\circ,V)$ is of $k$-eigenvalue type.
\end{theorem}

\begin{proof}
As $X$ has good reduction over $K$, there exists a rational integer
$M$ such that $X$ is the generic fiber of a smooth proper scheme
$\X$ over $\O_K[1/M]$, where $\O_K$ is the ring of integers of $K$.
Thus, the homomorphism $\Gal(\bar K/K)\to\Aut(H^k(\bar X,\Q_p))$
factors through $\rho\colon \Gamma_{K,Mp}\to \GL_n(\Q_p)$, the
Galois group over $K$ of the maximal subfield of $\bar K$ unramified
over any prime of $\O_K$ not dividing $Mp$.

For each prime $v$ of $\O_K$ dividing $Mp$, we fix an embedding $\bar
K\hookrightarrow \bar K_v$ and therefore an embedding
$\Gamma_{G_v}\hookrightarrow \Gamma_{K,Mp}$.  Let $G$, regarded as
an algebraic group over $\Q_p$, be the Zariski-closure of
$\rho(\Gamma_{K,Mp})$ in $\GL_n$, $G_v$ the Zariski-closure of
$\rho(\Gamma_{G_v})$, and $G_p$ the normal subgroup of $G$ generated
by $G_v^\circ$ for all $v$ lying over $p$. Replacing $K$ by a finite
extension, we may assume that $G_v$ is connected for all such $v$,
so $G_p$ is generated by conjugates of the $G_v$. By
Theorem~\ref{t:HT}, the complexification $G_{p\C}$, together with
its natural $n$ dimensional representation, is of $k$-eigenvalue
type.   If $G_p$ is of finite index in $G$, the theorem follows.
Otherwise, there exists a homomorphism $\Gamma_{K,Mp}\to
G(\Q_p)/G_p(\Q_p)$ with Zariski-dense, and therefore infinite
$p$-adic analytic image.  By construction, this homomorphism is
unramified at all primes over $v$.  Such a homomorphism cannot exist
according to the Fontaine-Mazur conjecture.

\end{proof}

\begin{cor}
If the Fontaine-Mazur conjecture is true, then for every complex
non-singular variety $X$ over a number field $K$, the Zariski
closure of the image of $\Gal(\bar K/K)$ in $\Aut(H^2(\bar X,\Q_p))$
has no factor of type $E_8$.
\end{cor}

\section{Application to braid group representations}

\label{s:appl} Artin's braid group $\B_m$ is generated by
$\sigma_1,\ldots,\sigma_{m-1}$ subject to relations
$$\sigma_i\sigma_j=\sigma_j\sigma_i \quad \text{if $|i-j|\geq 2$},\quad
\sigma_i\sigma_{i+1}\sigma_i=\sigma_{i+1}\sigma_i\sigma_{i+1} \text{
for $1\leq i\leq m-1$}.$$
 In \cite{FLW}, the closed images of the
unitary $q=e^{2\pi i/\ell}$ Hecke algebra representations of the
braid groups are completely analyzed (completing a program initiated
by Jones) for $\ell \ge 5$ and $\ell\neq 6$. In this section, we
will carry out a similar analysis.  We also discuss the situations
in which the braid group representations arising from quantum groups
at roots of unity satisfy the 3-eigenvalue condition.

\subsection{Set-up}\label{setup}
Given an irreducible unitary representation
 $(\rho,V)$ of $\B_m$ there
 are three distinct possibilities for $G=\overline{\rho(\B_m)}$
 \begin{enumerate}
 \item  $G/Z(G)$ is finite
 \item $\SU(V)\subset G$
 \item $G/Z(G)$ is infinite, but $\SU(V)\not\subset G$.
 \end{enumerate}
While the first (finite group) and third (non-dense) possibilities
are interesting, we will focus on the second.  There are a number of
reasons for doing this.  Firstly, we will see that $\SU(V)\subset G$
is the generic situation, while the other (non-dense) cases require
a case-by-case analysis that we will carry out in a separate work.
Also, density is crucial for
applications to quantum computing---our original motivation. Lastly,
the application of Theorem \ref{mainthm} leads most directly to the
conclusion $\SU(V)\subset G$, \emph{i.e.} by showing that $(G,V)$ is
an indecomposable pair satisfying the 3-eigenvalue property for
which the three eigenvalues do not form a geometric progression.
Nearly all of the finite group/non-dense examples come from pairs
having eigenvalues in geometric progression which will be considered
in a forthcoming paper by the first two authors.  We proceed with
the following program:
\begin{enumerate}
\item Determine which representations have exactly three
eigenvalues.
\item Determine conditions for the representations from (1) to be unitary.
\item Determine when the three eigenvalues from (1) and (2) satisfy the no-cycle
condition.  This will give us all pairs $(G,V)$.
\item Determine when the three eigenvalues from (1) and (2) are not in geometric
progression.  Although this does not ensure density, it does
guarantee the pair $(G,V)$ is indecomposable, as three eigenvalues
coming from a decomposable pair must be in geometric progression by
Lemma \ref{decomplemma}.
\item Determine when $G$ is infinite modulo the center for the cases not excluded by
(1)-(4).
\end{enumerate}

\subsection{$BMW$-algebra representations of the braid groups}
We apply the strategy outlined above to $BMW$-algebras, first
recalling what is well-known and then proceeding to the subsequent
steps.

\subsubsection{Definitions and combinatorial results}
Most of the material here can be found in \cite{We}, and we
summarize the details germane to the problem, carrying out steps (1)
and (2) in the above program.

The Birman-Wenzl-Murakami (BMW) algebras are a sequence of finite
dimensional algebras equipped with Markov traces.  They can be
described as quotients of the group algebra $\C(r,q)\B_m$ of Artin's
braid group where $r$ and $q$ are complex parameters.  The precise
definition of the BMW-algebra $\mathcal{C}_m(r,q)$ is:
\begin{defn}
Let $g_1, g_2,\ldots, g_{m-1}$ be invertible generators satisfying
the braid relations $(B1)$ and $(B2)$ above and:
\begin{enumerate}
\item[(R1)] $(g_i-r^{-1})(g_i-q)(g_i+q^{-1})=0$
\item[(R2)] $e_ig_{i-1}^{\pm 1}e_i=r^{\pm 1}e_i$, where
\item[(E)] $(q-q^{-1})(1-e_i)=g_i-g_i^{-1}$ defines $e_i$.
\end{enumerate}
\end{defn}

The relations (R2) can be best understood by pictures where $g_i$ is
the the braid generator $\sigma_i$ and $e_i$ is the $i$-th generator
of the Temperley-Lieb algebra.  Relation (R1) shows that the image
of $g_i$ in any representation of $\CC_m(r,q)$ has 3 eigenvalues:
$r^{-1},q$ and $-q^{-1}$. When $r\not=\pm q^n$ and $q$ is not a root of unity,
each $BMW$-algebra $\mathcal{C}_m(r,q)$ is finite-dimensional and
semisimple with simple components labeled by Young diagrams
with $m-2j\geq 0$ boxes for $j\in\N$. In other words, the
$BMW$-algebra is a direct sum of full matrix algebras. For each
simple component $\mathcal{C}_{m,\la}$ let $V_{m,\la}$ be the unique
non-trivial simple $\mathcal{C}_{m,\la}$-module. Then the branching
rule for restricting $V_{m,\la}$ to $\mathcal{C}_{m-1}(r,q)$ is:
$$V_{m,\la}\cong\bigoplus_{\mu\leftrightarrow\la}V_{m-1,\mu}$$
where $V_{m-1,\mu}$ is a simple $\mathcal{C}_{m-1}(r,q)$-module
 and $\mu$ is a Young diagram with $m-1-2j\geq 0$ boxes obtained
from $\la$ by adding/removing a box to/from $\la$.  This description
of inclusions among to $BMW$-algebras can be neatly encoded in a
graph called the \emph{Bratteli diagram}.
 The graph consists of vertices labelled by $(m,\la)$ with $|\la|=m-2k$ arranged in
 rows (labelled by integers $m$).
 Vertices in adjacent rows are connected if the their labels
 differ by 1 in the first entry and by one box in the second.  The dimension of $V_{m,\la}$
 can thus be computed by adding up the dimensions of the $V_{m-1,\mu}$ whose labels are
 connected to $(m,\la)$ by an edge. We obtain representations
of $\B_m$ on $\bigoplus_\la V_{m,\la}$ via the map
$\sigma_i\rightarrow g_i\in \CC_m(r,q)$.

We are interested in obtaining unitary representations of $\B_m$
from $BMW$-algebras, so we must consider semisimple quotients with
$r$ and $q$ specialized at roots of unity.  Specifically, we let
$r=q^n$ for $-1\not=n\in\Z$ and  $q=e^{\pi i/\ell}$ ($\ell\not=
1$), \emph{i.e.}, a primitive $2\ell$th root of unity.
 If a given irreducible representation is unitary for $q=e^{\pi i/\ell}$, it will remain so for
 $q=e^{-\pi i/\ell}$.  For other choices of primitive roots of unity
 we cannot expect to have unitarity.  The quotient of
each specialized $BMW$-algebra by the annihilator of the trace
$A_m:=\{a\in\CC_m(r,q): tr(ab)=0\quad \text{for all $b$}\}$ is semisimple
and we denote it $\CC_m(q^n,q)$ (where $q$ is understood to be
$e^{\pi i/\ell}$). The branching rules and simple decomposition
described above for the generic case still essentially apply to
$\CC_m(q^n,q)$, except that
 some components no longer appear, and fewer Young diagrams are needed
 to describe the persisting components (for all $m$).  Precisely which components
 survive depends on the values $\ell$ and $n$, and the derivation can be found in
 \cite{We}, the results of which we will describe below. For now it is enough to
 note that each simple component (sector) that does survive the quotient gives
 us an irreducible representation of $\B_m$.  Let $\rho^{(n,\ell)}_{(m,\la)}$
 acting on $V_{(m,\la)}^{(n,\ell)}$ be the representation of $\B_m$ corresponding to
 the simple component of $\CC_m(q^n,q)$ labeled by $\la$.  Since the conjugacy class
 of $\rho^{(n,\ell)}_{(m,\la)}(\sigma_1)$ generates the closed image of $\B_m$
 topologically, there is a chance that the pair
  $$(\overline{\rho_{(m,\la)}^{(n,\ell)}(\B_m)},V_{(m,\la)}^{(n,\ell)})$$ satisfies
  the 3-eigenvalue property.

As a first step we need to know the conditions under which the image
of $\sigma_1\in\B_m$ under $\rho_{(m,\la)}^{(n,\ell)}$ has 3
distinct eigenvalues. The answer is well-known to experts (see
\cite{We}): \emph{for $m\geq 3$, the image of $\sigma_1$ under the
irreducible representation $\rho_{(m,\la)}^{(n,\ell)}$ has 3
distinct eigenvalues precisely when $|\la|<m$ and $\CC_{3,\Box}$ is
three dimensional.} This is equivalent to the requirement that the
corresponding simple component $\CC_{m,\la}$ contains the simple
component $\mathcal{C}_{3,\Box}$. This is most easily seen by
considering the Bratteli diagram as described above.  It is shown in
\cite{We} that $\mathcal{C}_m(q^n,q)/\A_m\cong \I_m\oplus
\overline{\mathcal{H}}_m(q^2)$ where $\overline{\mathcal{H}}_m(q^2)$
is a quotient of the Iwahori-Hecke algebra of type $A_{m-1}$, and
$\I_m$ is the ideal generated by $e_{m-1}$ (see \cite{We}). The
Young diagrams labeling simple components of
$\overline{\mathcal{H}}_m(q^2)$ have $m$ boxes, whereas those of
$\I_m$ have $m-2j$ boxes for some $j\geq 1$. The representations of
$\B_m$ corresponding to the Iwahori-Hecke algebra part of
$\CC_m(q^n,q)$ have been studied in \cite{FLW} where they are
analyzed using the solution to the 2-eigenvalue problem.  Thus the
image of $\sigma_1$ on the irreducible representation $V_{m,\la}$
($m\geq 3$) has (exactly) 3 distinct eigenvalues precisely when
$|\la|<m$ and $\CC_{3,\Box}$ is 3-dimensional in which case the
eigenvalues are $\{q^{-n},q,-q^{-1}\}$. We can eliminate many
redundant cases using isomorphisms (see \cite{TubaWenzl2}):
\begin{equation}\label{isos}
\CC_m(q^n,q)\cong\CC_m(-q^{-n},q)\cong\CC_m(-q^n,-q)\cong\CC_m(q^{-n},q^{-1}).
\end{equation}
We describe the restrictions more precisely in the following, which
is a reformulation of several results in \cite{We} and \cite{R}.
Denote by $\la_i$ (resp. $\la_i^\prime$) the number of boxes in the
$i$th row (resp. column)
 of the Young diagram $\la$.
\begin{prop}\label{wenzlprop}
 Let $q=e^{\pi i/\ell}$ and $m\geq 3$.
\begin{enumerate}
\item  The matrix algebra $\CC_{3,\Box}$ is a simple 3-dimensional
 subalgebra of $\CC_m(q^n,q)$, provided one of the following conditions holds:
\begin{enumerate}
\item $n=1$ and $\ell\geq 3$
\item $n=2$ and $\ell\geq 4$
\item $3\leq n\leq \ell-3$ (so $\ell\geq 6$)
\item $4-\ell\leq n\leq -4$, $n$ is even and $\ell$ is odd (so $\ell\geq 9$)
\item $5-\ell\leq n\leq -5$, $n$ is odd and $\ell$ is even (so $\ell\geq 10$)
\end{enumerate}
Moreover, this list is exhaustive up to the isomorphisms \ref{isos}.
\item The $\la$ for which $\CC_{m,\la}$ may appear as a simple component
in some $\CC_m(q^n,q)$ are in the following sets of
\textbf{$(n,\ell)$-admissible} Young diagrams
 corresponding to each of the 5 cases above:
\begin{enumerate}
\item $\{[1^2]\}\cup\{[k]: k\in\N\}$
\item $\{[1^3]\}\cup\{[k],[k,1]: 1\leq k\leq \ell-1\}$
\item $\{\la: \la_1+\la_2\leq\ell-n+1 \quad \text{and}
\quad \la_1^\prime+\la_2^\prime\leq n+1\}\cup\{[\ell-n+1,1^{n-1}]\}$
\item $\{\la: \la_1+\la_2\leq 1-n \quad \text{and} \quad \la_1^\prime\leq (\ell+n-1)/2\}$
\item $\{\la: \la_1\leq(-1-n)/2 \quad \text{and} \quad \la_1^\prime\leq (\ell+n-1)/2\}$
\end{enumerate}
\item Thus the image of $\sigma_1$ under the irreducible
representation $\rho^{(n,\ell)}_{(m,\la)}$ with $|\la|<m$ has 3
distinct eigenvalues provided $n$ and $\ell$ satisfy one of the
conditions of (1) and $\la$ is in the corresponding set of
admissible Young diagrams in (2).  These representations are unitary
except possibly in case (d).
\end{enumerate}
\end{prop}
\begin{remark}
Observe that the set in 2(a) is infinite and independent of $\ell$.
The other four labeling sets are finite, and it is easy to see that
the corresponding Bratteli diagrams are periodic. In the case $n=2$ there
is a slight exception to the rule for constructing the Bratteli diagram:
the diagrams labeled by $[\ell-1,1]$ and $[\ell-1]$ are \emph{not} connected
by an edge (see \cite{We}, Prop. 6.1). The fact that the
representations in (a),(b),(c) and (e) are unitary was proved in
\cite{We}.  The full (reducible) representations of $\B_m$ factoring
over $\CC_m(q^n,q)$ corresponding to case (d) were shown in \cite{R}
to be non-unitarizable \emph{for any $q$} when $\ell> 2(-n+1)$. This
leaves only finitely many possible $\ell$ for each fixed $n$, and
even in these cases one can use the techniques of \cite{R} to show
that for $q=e^{\pi i/\ell}$ one does not get unitarity except in
degenerate cases.  Restricting to the irreducible sectors one may
get unitarizable representations, but not uniformly, so that for
$m\gg 0$ no irreducible sector is unitary.
\end{remark}

\subsubsection{Cycles and geometric progressions} The eigenvalues of
any of the irreducible representations satisfying the conditions of
Proposition \ref{wenzlprop} are $\{q,-q^{-1},q^{-n}\}$, with
$q=e^{\pi i/\ell}$.  Steps (3) and (4) of the program can be
accomplished with simple computations.  We have:
\begin{lemma}\label{nocyclegp}
Let $n$, $\ell$  and $\la$ be as in Proposition \ref{wenzlprop}.
Then the eigenvalues of $\rho_{(m,\la)}(\sigma_1)$:
\begin{enumerate}
\item
fail the no-cycle property if and only if $n=1$ or $(n,\ell)=(3,6)$
and
\item are in geometric progression if and only if
$n\in\{3,\ell-3,\pm\ell/2\}$.
\end{enumerate}
\end{lemma}
\begin{proof} The only way $\{q,-q^{-1},q^{-n}\}$ can fail the no-cycle condition
is if it contains a coset of $\{\pm 1\}$ or $\{1,e^{2\pi
i/3},e^{4\pi i/3}\}$.  With the restrictions in Prop.
\ref{wenzlprop} that $\ell\geq n+2$ for $n>0$ and $\ell\geq 4-n$ for
$n<0$ as well as $\ell\geq 3$ one checks that only $n=1$ and
$(n,\ell)=(3,6)$ fail no-cycle.  For the eigenvalues to be in
geometric progression (still satisfying the conditions of Prop.
\ref{wenzlprop}) we check the solutions of $\la_1\la_2-(\la_3)^2=0$
for the three possible assignments of $\la_3$.  These yield the
three solutions for $n$ above. \end{proof}

\begin{remark}
All of the exceptional cases $n\in\{1,3,\ell-3,\pm\ell/2\}$ will be
considered in a future work.  As we remarked above the case $n=1$ is
unique in that the labelling set of irreducible sectors is infinite.
 In fact, it is not hard to see, using the classification of $m$-dimensional
 irreducible representations of $\B_m$ found in \cite{formanek}, that
 one obtains some finite group images for every $m$ when $n=1$.
 By the isomorphisms of $BMW$-algebras corresponding to $r\leftrightarrow -r^{-1}$
 we see that the two cases $n=3$ and $n=\ell-3$ are actually the
 same.  Moreover, it can be shown that the (specialized quotient)
 $BMW$-algebras $\CC_m(q^3,q)$ can be embedded (diagonally) in quotients of the
tensor squares of Iwahori-Hecke algebras $\mathcal{H}_m(q^2)$.  This
indicates that the corresponding pairs may be tensor decomposable.
In the subcase $(n,\ell)=(3,6)$ work of Jones in \cite{J1} shows
that the images are all finite groups (essentially $\PSL(2m,3)$). The
case $n=-\ell/2$ sometimes also have finite group images \emph{e.g.}
when $(n,\ell)=(-5,10)$, see \cite{J2}.
\end{remark}

\subsubsection{Infinite images and density}
 Finally, we need to determine, for representations not excluded by the
 steps (1)-(4) above, the values of $m$, $\ell$, $n$, and
$\la$ for which the image of $\B_m$ under the unitary irreducible
representation $\rho^{(n,\ell)}_{(m,\la)}$ in $\CC_m(q^n,q)$ with
$q=e^{\pi i/\ell}$ is infinite modulo the center.  Proposition
\ref{wenzlprop} implies that a sufficient condition for
$\rho^{(n,\ell)}_{(m,\la)}$ to have infinite image is that the
3-dimensional representation $\rho^{(n,\ell)}_{(3,\Box)}$ have
infinite image.  So as a first step, we study this condition. For
convenience of notation we denote this representation simply by
$\rho$ despite its dependence on the parameters. A non-unitary
realization of $\rho$ is given by:
$$\sigma_1 \rightarrow A:=\begin{pmatrix}\frac{1}{q^{n}} & \frac{q^2-1}{q} &
0\\0 & \frac{q^2-1}{q}& i\\
0 & -i & 0\end{pmatrix}, \sigma_2\rightarrow B:=\begin{pmatrix}0 &0&
-i\\
0 &\frac{1}{q^n}& \frac{-i(q^2-1)}{q^{n+1}}\\ i & 0&
\frac{q^2-1}{q}\end{pmatrix}$$ found in \cite{BW}.

  Blichfeldt \cite{Bl} has determined the irreducible finite
  subgroups of $\PSL(3,\C)$.  Six are primitive groups of orders 36, 60,
72, 168, 216, and 360, and the imprimitive subgroups come in two
infinite families isomorphic to extensions of $S_3$
and $\Z_3$ by abelian groups.
\begin{defn}
A group $\Gamma$ is \emph{primitive} if $\Gamma$ has a faithful irreducible
representation which cannot be expressed as a direct sum of
subspaces which $\Gamma$ permutes nontrivially.
\end{defn}
By Lemma \ref{l:permute}, a sufficient condition for $G=\overline{\rho(\B_3)}$ to
be primitive is that the spectrum of $\rho(\sigma_1)$ satisfies the
no-cycle property.  So by Lemma \ref{nocyclegp} the image of $\rho$
is only imprimitive in the excluded cases $n=1$ and
$(n,\ell)=(3,6)$.
So we may assume that the $G$ is primitive.  We wish
to determine when $G$ is infinite modulo the center.  By rescaling the
images of the generators $\sigma_i$ by the cube root of the determinant of $\rho(\sigma_i)$
we may assume that $G\subset \SL(3,\C)$, and
to determine the image modulo the center it suffices to consider the projective image.
Thus $G/Z(G)\subset \PSL(3,\C)$, and we may apply Blichfeldt's classification.  We state his result
and include some useful information about orders of elements in:
\begin{prop}
The primitive subgroups of $\PSL(3,\C)$ are:
\begin{enumerate}
\item The \emph{Hessian} group $H$ of order 216 or a normal subgroup of $H$ of order 36 or 72.
The Hessian group is the subgroup of $A_9$ generated by $(124)(568)(397)$ and $(456)(798)$, and has
elements of order $\{1,2,3,4,6\}$.
\item The simple group $\PSL(2,7)\subset A_7$ of order 168.  The orders of elements are $\{1,2,3,4,7\}$.
\item The simple group $A_5$ having elements of orders $\{1,2,3,5\}$.
\item The simple group $A_6$ having elements of orders $\{1,2,3,4,5\}$.
\end{enumerate}
\end{prop}

Using this result we have the following:
\begin{theorem}\label{BMWthm}

Let $n$ and $\ell$ be chosen so that $\rho^{(n,\ell)}_{(3,\Box)}$ is a 3-dimensional unitary irreducible
representation of $\B_3$ with eigenvalues not in geometric progression and satisfying the no-cycle condition.
That is, $n$ and $\ell$ satisfy the hypotheses of Proposition \ref{wenzlprop}(1)(b),(c) or (e) in addition to
$n\not\in\{3,\ell-3,\pm\ell/2\}$.  Let $m\geq 3$ and $|\la|<m$ with $\la$
$(n,\ell)$-admissible.
 The closure of the group $\rho^{(n,\ell)}_{(m,\la)}(\B_m)$ is infinite
 modulo the center with two exceptions: if
 $(n,\ell)\in\{(-5,14),(-9,14)\}$ with $(m,\la)\in\{(3,\Box),(4,[0])\}$
 then the projective images are isomorphic to $\PSL(2,7)$.
Excluding these cases, if the dimension of the representation $\rho^{(n,\ell)}_{(m,\la)}$ is $k$,
then the closure of the image of $\B_m$ contains $SU(k)$.
 \end{theorem}
\begin{proof} Knowing the specific eigenvalues of $\rho(\sigma_1)$
we compute its projective order $t(n,\ell)$ as a function of $\ell$
and $n$ to be:
\begin{equation}
t(n,\ell)=
\begin{cases} \ell/2 & \text{if $\ell\equiv 2\pmod{4}$ and $n\equiv 3\pmod{4}$}\\
\ell & \text{if $\ell\equiv 0\pmod{4}$ and $n$ even or}\\ & \text{
$\ell\equiv 2\pmod{4}$ and
$n\equiv 1\pmod{4}$}\\
2\ell &\text{otherwise}\end{cases}
\end{equation}
Under the stated hypotheses on $n$ and $\ell$ we consider cases, comparing with
the list of possible orders of elements in Blichfeldt's classification.
\begin{enumerate}
\item If $\ell$ is odd, then $\ell\geq 5$ in which case $t(n,\ell)\geq 10$ which is too large.
\item If $\ell\equiv 0\pmod{4}$ then $\ell=8$ is the smallest value not yet excluded which gives
$t(n,\ell)\geq 8$ which is again too large.
\item If $\ell\equiv 2\pmod{4}$ then $\ell\geq 6$ and $t(n,\ell)\geq
12$ unless $n$ is odd.  If $n\equiv 1\pmod{4}$ then $\ell\geq 10$
which gives us $t(n,\ell)=\ell\geq 10$ which does not appear on the
list.  When $\ell\equiv 2\pmod{4}$ and $n\equiv 3\pmod{4}$ with $n>0$
we must have $n\geq 7$ which forces $\ell\geq 18$ since
$n\not=\ell/2$. For $\ell\equiv 2\pmod{4}$ and $n\equiv 3\pmod{4}$
with $n<0$ we must have $\ell\geq 14$ (since $\ell=10$ leads to
$n=-5=-\ell/2$), which has the two possible values $n=-5$ or $n=-9$
which we claim gives rise to finite images.  Observe that
$t(-5,14)=t(-9,14)=7$.
\end{enumerate}
To show that the projective images for $(-5,14)$ and $(-9,14)$ are
both $\PSL(2,7)$ we first observe that it is enough by the
isomorphism of \ref{isos} with $r\leftrightarrow -r^{-1}$ and
$q^{-5}\leftrightarrow q^{-14+5}=q^{-9}$ so these two cases give the
same images.  Then we use the explicit matrices $A$ and $B$ above to
define $S=B^{-1}$ and $T=BAB$ which then (projectively) satisfy the
relations $S^7=(S^4T)^4=(ST)^3=T^2=I_{3\times 3}$ defining
$\PSL(2,7)$.  It is immediate from the Bratteli diagram that the
representation of $\B_4$ corresponding to $(4,[0])$ is irreducible and isomorphic
to the that of $(3,\Box)$ when
restricted to $\B_3$.  Moreover, the representations of
$\B_4$ corresponding to diagrams $[1^2]$ and $[2]$ each contain the
representation of $\B_3$ corresponding to the Young diagram
$[1^2,1]$ which was shown in \cite{FLW} to have infinite image
(modulo the center).  For all of the infinite image cases the hypotheses
of Theorem \ref{mainthm} are satisfied and the eigenvalues are not in
geometric progression so density follows.
\end{proof}

\subsection{Quantum Groups}
In this subsection we consider braid group actions on centralizer
algebras of representations of quantum groups at roots of unity.  We
find and analyze examples in which
 we may apply Theorem \ref{mainthm}.  We follow the general strategy in
 Subsection \ref{setup}, but we note that as the
 representation spaces available to us are not necessarily simple subquotients
 of braid group algebras (unlike $BMW$-algebras) there is a subtlety
 regarding irreducibility.

\subsubsection{Braid group action on centralizer algebras}
The Drinfeld-Jimbo quantum group $U:=U_q\g$ associated to a simple
Lie algebra $\g$ is a ribbon Hopf-algebra.  The so-called
\emph{universal $R$-matrix} that intertwines the coproduct with the
opposite coproduct on $U$ can be used to construct representations
of the braid group $\B_n$ on the morphism space $\End_U(\Vn)$ for
any finite dimensional highest weight $U$-module $V$ as follows.
Fix such a $U$-module $V$ and define $\check{R}=P_V\circ
R\mid_{V\otimes V}\in\End_U(V^{\otimes 2})$ to be the
$U$-isomorphism afforded us by composing the image of the universal
$R$-matrix acting on $V\otimes V$ with the ``flip" operator
$P_V:v_1\otimes v_2\rightarrow v_2\otimes v_1$.  Then define
isomorphisms for each $1\leq i\leq n-1$ :
$$\check{R}_i:=\one^{\otimes
(i-1)}\otimes\check{R}\otimes\one^{\otimes (n-i-1)}\in\End_U(\Vn)$$
so that the $\check{R}_i$ satisfy the braid group relations. Then
define a representation of $\B_n$ on $\End_U(\Vn)$ by
$\sigma_i.f=\check{R}_i\circ f$.

Lusztig has defined a modified form
of the quantum group $U$ so that one may specialize the quantum
parameter $q$ to $e^{\pm\pi i/\ell}$.  In fact, one may choose any $q$
so that $q^2$ is a primitive $\ell$th root of unity, but we will
restrict our attention to $q=e^{\pi i/\ell}$ since these values
(sometimes) yield unitary representations (see \cite{wenzlcstar}),
 which remain unitary for
$\overline{q}$.  The full representation category of $U$ at roots of
unity is not semisimple, but has a semisimple subquotient category.
This process is essentially due to Andersen and his coauthors (see
\cite{An} and references therein).  This yields a semisimple ribbon
category $\F$ (see \cite{Turaev} for the definitions) with finitely
many simple objects labeled by highest weights in a truncation of
the dominant Weyl chamber, called the \emph{Weyl alcove}.  The braid
group still acts on $\End_U(\Vn)$ for any object $V$ as above, and
for each quantum group we look for simple objects $V_\la$ so that
the images of the braid generators on the irreducible
subrepresentations of $\End_U(V_\la^{\otimes n})$ have 3
eigenvalues.  Because the tensor product rules for objects labelled
by weights near the upper wall of the Weyl alcove depends on $\ell$,
we do not explicitly determine all $V_\la$ giving rise to pairs with
the 3-eigenvalue property, and restrict our attention to weights
near $0$.   As in the $BMW$ algebra setting, we will always have an
irreducible 3-dimensional representation of $\B_3$ to which we may
reduce most questions. We sketch the idea (see \emph{e.g.}
\cite{TubaWenzl1} Section 3): If $V$ is a simple object in (a finite
semisimple ribbon category) $\F$ such that $V\otimes V\cong
V_0\oplus V_1\oplus V_2$ is the decomposition into 3 inequivalent
simple objects then $\End_U(V^{\otimes 3})$ has a 3-dimensional
irreducible subrepresentation isomorphic to $\Hom_U(V^{\otimes
3},W)$ for a simple object $W$ appearing in $V^{\otimes 3}$ with
multiplicity three.  Provided the (categorical) $q$-dimension of
each of $W$, $V$ and $V_i$ are non-zero then this representation is
irreducible and the image of $\sigma_1$ will have three distinct
eigenvalues. As in the $BMW$-algebra situation we can construct a
Bratteli diagram encoding the containments of the semisimple finite
dimensional algebras:
$$\End_U(V)\subset\End_U(V\otimes V)\subset\cdots\subset\End_U(\Vn)\cdots$$
 The simple components of $\End_U(\Vn)$ will be isomorphic
to $\Hom_U(\Vn,V_\mu)$ where $V_\mu$ is a simple object appearing in
the decomposition of $\Vn$.  The edges of the Bratteli diagram are
determined by decomposing $V_\gamma\otimes V$ where $V_\gamma$ is a
simple subobject of $V^{\otimes (n-1)}$.  There are techniques known
for obtaining these decompositions, for example Littelman's path
basis technique \cite{Li}, or crystal bases. However, when we
consider the action of the braid group $\B_n$ on the spaces
$\Hom_U(\Vn,V_\mu)$ we have no guarantee that the action is
irreducible.  This is because $\End_U(\Vn)$ might not be generated
by the image of $\B_n$.

\subsubsection{Density results}
We proceed to find pairs $(X_r,\la)$ so
that the ribbon category corresponding to the quantum group
$U_q\g(X_r)$ of Lie type $X_r$ has simple object $V_\la$ with
$V_\la^{\otimes 3}\cong V_0\oplus V_1\oplus V_2$ as above.  We find
that $(A_r,\varpi_2)$, $(A_r,2\varpi_1)$, $(B_r,\varpi_1)$,
$(C_r,\varpi_1)$, $(D_r,\varpi_1)$ and $(E_6,\varpi_1)$ do satisfy
these conditions (where the weights $\varpi_i$ are labeled as in
\cite{Bo}). With these in hand we compute the eigenvalues of the
images of $\sigma_i$ in the corresponding representations. We use
the following result found in \cite{LeducRam}, Corollary 2.22(3)
originally due to Reshetikhin.  The form $\lan\cdot,\cdot\ra$
is the symmetric inner product on the root lattice normalized so
that the square lengths of \emph{short} roots is 2, and the
 weight $\rho$ is the half
sum of the positive roots.
\begin{prop}\label{qeigs}
Suppose that $V=V_\varpi$ is an irreducible representation of the
quantum group $U_q\g$ and that $V\otimes V_\la$ is multiplicity free
for all $V_\la$ appearing in some $\Vn$.  Then for any $V_\nu$
appearing in $V^{\otimes 2}$ we have:
$$\check{R}_i\mid_{V_\nu}=\pm q^{(1/2)\lan
\nu,\nu+2\rho\ra-\lan\varpi,\varpi+2\rho\ra}\one_{V_\nu}$$ where the
sign is $+1$ if $V_\nu$ appears in the symmetrization of $V^{\otimes
2}$ and $-1$ if $V_\nu$ appears in the anti-symmetrization of
$V^{\otimes 2}$.
\end{prop}
We record the results in Table 1, where the notation follows
\cite{Bo}.  The symbol $\1$ denotes the unit object in the category.
The necessary computations are standard and can be done by hand
\emph{e.g.} using the technique of \cite{Li}.
\begin{table}\label{qgroup}
\begin{tabular}{*{4}{|c}|}
\hline
$(X_r,\la)$ & $S^2(V_\la)$ & $\bigwedge^2(V_\la)$ & Eigenvalues\\
\hline\hline $(A_r,\varpi_2)$ & $V_{2\varpi_2}$ &
$V_{\varpi_1+\varpi_3}\oplus V_{\varpi_4}$&
$q^{\frac{4}{r+1}+1}\{q,-q^{-1},-q^{-5}\}$\\
\hline
$(A_r,2\varpi_1)$ & $V_{2\varpi_2}\oplus V_{4\varpi_1}$ & $V_{2\varpi_1+\varpi_2}$
& $-q^{\frac{4}{r+1}-1}\{-q^{-1},-q^{5},q\}$ \\
\hline
$(B_r,\varpi_1)$ & $V_{2\varpi_1}\oplus\1$ & $V_{\varpi_2}$ & $\{q^2,q^{-4r},-q^{-2}\}$\\
\hline
$(C_r,\varpi_1)$ &$V_{2\varpi_1}$ &$V_{\varpi_2}\oplus\1$& $\{q,-q^{-1},-q^{-2r-1}\}$ \\
\hline
$(D_r,\varpi_1)$ & $V_{2\varpi_1}\oplus\1$ & $V_{\varpi_2}$ & $\{q, q^{2r-1},-q^{-1}\}$\\
\hline
$(E_6,\varpi_1)$ & $V_{2\varpi_1}\oplus V_{\varpi_6}$& $V_{\varpi_3}$  & $q^{1/3}\{q, q^{-9},-q^{-1}\}$\\
\hline
\end{tabular}\caption{Eigenvalues of $\check{R}_i$}
\end{table}
 The braid group representations corresponding to Lie types $B, C$ and $D$ are the
 same as those factoring over specializations of $BMW$-algebras, due to $q$-Schur-Weyl-Brauer
 duality, see \cite{We}.  For this reason we ignore these cases in the following
 weaker version of Theorem \ref{BMWthm}.
\begin{theorem}
Let $(X_r,\la)$ be as in Table 1 with $X=A_r$ or $E_6$.  Then
\begin{enumerate}
\item For $(A_r,\varpi_2)$:  provided $r\geq 3$ and $\ell\geq \rm{max}(r+3,7)$,
$\Hom_U((V_\la)^{\otimes 3},V_{\varpi_2+\varpi_4})$ is unitary, irreducible
and 3-dimensional and the image of $\sigma_1$ has 3 distinct eigenvalues.  If $V_\mu$ appears in
$V_{\varpi_2+\varpi_4}\otimes V_\la^{\otimes n-3}$ then $\Hom_U(\Vn,V_\mu)$ \emph{contains}
an irreducible unitary representation of $\B_n$
with the 3-eigenvalue property.  When $\ell\not\in\{10,14\}$,
the eigenvalues of the image of $\sigma_1$ are not in geometric progression and the images of $\B_n$
are infinite modulo the center and so are dense in these cases.
\item For $(A_r,2\varpi_1)$: $\Hom_U((V_\la)^{\otimes 3},V_{2\varpi_1+2\varpi_2})$ is unitary, irreducible
and 3-dimensional provided $r\geq 1$ and $\ell\geq r+5$. If $V_\mu$
appears in $V_{2\varpi_1+2\varpi_2}\otimes V_\la^{\otimes n-3}$ then
$\Hom_U(\Vn,V_\mu)$ \emph{contains} an irreducible unitary
representation of $\B_n$ with the 3-eigenvalue property. When
$\ell\not\in\{6,10\}$ the eigenvalues of the image of $\sigma_1$ are
not in geometric progression and the images of $\B_n$ are infinite
modulo the center and so are dense in these cases.
\item for $(E_6,\varpi_1)$: $\Hom_U((V_\la)^{\otimes 3},V_{\varpi_1+\varpi_6})$ is unitary, irreducible
and 3-dimensional provided $\ell\geq 14$.  If $V_\mu$ appears in
$V_{\varpi_1+\varpi_6}\otimes V_\la^{\otimes n-3}$ then
$\Hom_U(\Vn,V_\mu)$ \emph{contains} an irreducible unitary
representation of $\B_n$ with the 3-eigenvalue property.  Provided
$\ell\not=18$, the eigenvalues of the image of $\sigma_1$ are not in
geometric progression and the images of $\B_n$ are infinite modulo
the center and so are dense in these cases.
\end{enumerate}
\end{theorem}
\begin{proof} For the object labelled by $V_\nu$ to be in the fundamental alcove, we must have
$\lan \nu+\rho, \theta\ra<\ell$ where $\theta$ is the highest root.
This condition together with the requirement that the eigenvalues be
distinct yield the first restrictions in each case.  The unitarity
of the representations is shown in \cite{wenzlcstar}.  In each case
the representation spaces $\Hom_U(\Vn,V_\mu)$ described in the
theorem contain the 3-dimensional representation spaces, so by
restriction to $\B_3$ we see that the $\B_n$ representations must
contain an irreducible unitary subrepresentation with the
3-eigenvalue property.  Geometric progressions appear in each of the
three cases if and only if $\ell=10$ in the first case, $\ell=6$ or
$10$ in the second case and $\ell=18$ in the last case.  Computing
projective orders of the images of $\sigma_1$ and comparing as in
the proof of Theorem \ref{BMWthm} we find that the only finite group
image that arises is in the first case with $\ell=14$.  With these
exceptions, the hypotheses of Theorem \ref{mainthm} are satisfied
and we may conclude the images are dense.
\end{proof}
\begin{remark}
To get sharper results we would need to describe the decompositions
of the $\B_n$ representations $\Hom_U(\Vn,V_\mu)$ that appear in the
above theorem.  This is in general quite complicated. In fact, the
type $E_6$ case appears in an exceptional series discussed in
\cite{WenzlE} (and extended slightly in \cite{R2}).  These give new
semisimple finite dimensional quotients of the braid group algebras
analogous to $BMW$-algebras about which little is known.
\end{remark}

\subsection{Concluding Remarks}

In comparing this work to the 2-eigenvalue paper,
it may be noted that we do not provide applications of our results
to the distribution of values of the Kauffman polynomial
in analogy with those given for the Jones polynomial in  \cite[\S5]{FLW}.
That is, we do not consider the set of values $F_L(a,z)$ for fixed $a$ and $z$
and varying $L$.  These values can be described as linear combinations of
traces of any braid with closure $L$ in the different irreducible factors of a
BMW algebra, just as in \cite{FLW}.
The difficulty is that our information on the closures of
braid groups in BMW algebras is less detailed than the corresponding
information for Hecke algebras.  In particular, we have not
completely determined these closures for the irreducible factors
of the BMW-representations which are excluded in the statement of
Theorem~\ref{BMWthm}.  Neither have we determined the equivalences
and dualities existing between different irreducible factors in a fixed
BMW-algebra.  We certainly expect the limiting distributions to be Gaussian
as for the Jones polynomial,
but we do not yet have enough information to ensure that this is so.

In 1990s, Vertigan (see Theorems 6.3.5 and 6.3.6 of \cite{Wel})
analyzed the classical computational complexity of exactly evaluating
various knot polynomials at fixed complex values.
With a few exceptions, all evaluations are $\# P$-hard.  At these
few exceptional values, the link invariants have classical
topological interpretations and can be computed in polynomial time.
These results fit very well with the analysis of closed images of the braid
group representations. In the case of unitary Jones representations
of the braid groups at $q$, the closed image is dense in the
corresponding special unitary groups exactly when computing the link
invariants is $\#P$-hard at $q$, while the finite
image cases correspond to polynomial time computations.
Part of the appeal of working out the exceptions to Theorem~\ref{BMWthm}
is the hope of relating these cases to interesting special values of
the Kauffman polynomial.

\end{document}